\documentclass[12pt]{amsart}
\usepackage{amssymb,latexsym,comment,url}

\newcommand{\adj}{\operatorname{adj}}
\newcommand{\C}{{\mathbb C}}
\newcommand{\CC}{{\mathcal C}}
\newcommand{\Char}{\operatorname{char}}
\newcommand{\Diag}{\operatorname{Diag}}
\newcommand{\Div}{\operatorname{Div}}
\newcommand{\dv}{{\,\mid\,}}
\newcommand{\F}{{\mathbb F}}
\newcommand{\GL}{\operatorname{GL}}
\newcommand{\Jac}{\operatorname{Jac}}
\newcommand{\Kbar}{\overline{K}}
\newcommand{\la}{\lambda}
\newcommand{\level}{\ell}
\newcommand{\OK}{{\mathcal{O}_K}}
\newcommand{\Pf}{\operatorname{Pf}}
\newcommand{\PP}{{\mathbb P}}
\newcommand{\Q}{{\mathbb Q}}
\newcommand{\R}{{\mathbb R}}
\newcommand{\rank}{\operatorname{rank}}
\newcommand{\ra}{{\longrightarrow}}
\newcommand{\sh}{{\text{\rm sh}}}
\newcommand{\Sing}{\operatorname{Sing}}
\newcommand{\SL}{\operatorname{SL}}
\newcommand{\W}{{\mathcal W}}
\newcommand{\X}{{\mathcal X}}
\newcommand{\Z}{{\mathbb Z}}

\newenvironment{Proof}{\par\noindent{\sc Proof:}}%
                      {\hspace*{\fill}\nobreak$\Box$\par\medskip}
\newenvironment{ProofOf}[1]{\par\noindent{\sc Proof of #1:}}%
                       {\hspace*{\fill}\nobreak$\Box$\par\medskip}

\newtheorem{Proposition}{Proposition}[section]
\newtheorem{Theorem}[Proposition]{Theorem}
\newtheorem{Lemma}[Proposition]{Lemma}
\newtheorem{Corollary}[Proposition]{Corollary}

\theoremstyle{definition}
\newtheorem{Definition}[Proposition]{Definition}
\newtheorem{Remark}[Proposition]{Remark}
\newtheorem{Example}[Proposition]{Example}

\begin{document}

\title[Minimisation and reduction]%
{Minimisation and reduction of 5-coverings of elliptic curves}

\author{Tom~Fisher}
\address{University of Cambridge,
         DPMMS, Centre for Mathematical Sciences,
         Wilberforce Road, Cambridge CB3 0WB, UK}
\email{T.A.Fisher@dpmms.cam.ac.uk}

\date{21st December 2011}  

\begin{abstract}
We consider models for genus one curves of degree~$5$, 
which arise in explicit $5$-descent on
elliptic curves.  We prove a theorem on the existence of minimal models
with the same invariants as the minimal model of the Jacobian elliptic
curve and give an algorithm for computing such models.
Finally we describe how to reduce genus one
models of degree~$5$ defined over $\Q$.
\end{abstract}

\maketitle

\renewcommand{\baselinestretch}{1.1}
\renewcommand{\arraystretch}{1.3}
\renewcommand{\theenumi}{\roman{enumi}}

\section*{Introduction}
\label{intro}

Let $E$ be an elliptic curve defined over a number field~$K$.  An 
{\em $n$-covering} of~$E$ is a morphism $\pi : C \to E$ where $C$ is a
smooth curve of genus one, and $\pi = [n] \circ \psi$ for some
isomorphism $\psi : C \to E$ defined over $\Kbar$.  An $n$-descent on
$E$ computes the everywhere locally soluble $n$-coverings of $E$.  For
such $n$-coverings we have $\psi^*(n. 0_E) \sim D$ for some
$K$-rational divisor $D$ on $C$. The complete linear system $|D|$
defines a morphism $C \to \PP^{n-1}$. Thus in the cases $n = 2,3,4$ we
may represent $C$ by a binary quartic, ternary cubic, or pair of
quadrics in $4$ variables.  In the case $n=5$ we obtain curves
$C \subset \PP^4$ of degree $5$ that are defined by the $4 \times 4$ 
Pfaffians of a $5 \times 5$ alternating matrix of linear forms.

The question naturally arises as to how we can choose co-ordinates on
$\PP^{n-1}$ so that the equations for $C$ have small coefficients.
In the cases $n=2,3,4$ this was answered in \cite{CFS}, using the combination
of two techniques called {\em minimisation} and {\em reduction}. In this
paper we extend to the case $n=5$.  We establish results on minimisation
over an arbitrary local field (immediately implying results over 
any number field of class number~$1$), whereas those for reduction are 
specific to the case $K=\Q$. Implementations of our algorithms in the case 
$K=\Q$ are available in {\sf MAGMA} \cite{magma}. 

\section{Genus one models}
\label{g1mod}

A {\em genus one model} (of degree 5) is a $5 \times 5$ alternating
matrix of linear forms in variables $x_1, \ldots, x_5$.  We write $X_5(R)$ for
the space of all genus one models with coefficients in a ring $R$.
Models $\Phi$ and $\Phi'$ are {\em $R$-equivalent} if $\Phi' = [A,B]
\Phi$ for some $A,B \in \GL_5(R)$. Here the action of $A$ is via $\Phi
\mapsto A \Phi A^T$, and the action of $B$ is via $(\Phi_{ij}(x_1,
\ldots,x_5)) \mapsto (\Phi_{ij}(x'_1, \ldots,x'_5))$ where $x'_j =
\sum_{i=1}^5 B_{ij} x_i$.  The {\em determinant} of the transformation $g =
[A,B]$ is $\det g = (\det A)^2 \det B$.

We write $\Pf(\Phi)$ for the row vector $(p_1, \ldots, p_5)$ where
$p_i$ is $(-1)^{i-1}$ times the Pfaffian of the $4 \times 4$ submatrix
obtained by deleting the $i$th row and column of $\Phi$. This choice
of signs is made so that $\Pf(\Phi) \Phi = 0$. For $A \in \GL_5(R)$ we
note that $\Pf ( A \Phi A^T ) = \Pf(\Phi) \adj A$.

A genus one model $\Phi \in X_5(K)$ over a field $K$ is 
{\em non-singular} if the subscheme 
$\CC_{\Phi} = \{ \rank \Phi \le 2 \} \subset \PP^4$ defined
by the $4 \times 4$ Pfaffians of $\Phi$ is a smooth curve of genus
one.  We write $K[X_5]$ for the polynomial ring in the $50$ coefficients
of a genus one model. A polynomial $F \in K[X_5]$ is an {\em invariant}
of {\em weight} $k$ if $F \circ g = (\det g)^k F$ for all $g = [A,B]$
with $A,B \in \GL_5(\Kbar)$.  Taking $A$ and $B$ 
to be scalar matrices shows that an invariant of weight $k$ is a 
homogeneous polynomial of degree $5k$.

\begin{Theorem}
\label{thm:invjac}
Let $c_4, c_6, \Delta \in \Z[X_5]$ be the invariants of weights $4,6,12$,
satisfying $c_4^3 - c_6^2 = 1728 \Delta$, and scaled as specified in 
\cite{g1inv}. 
  \begin{enumerate}
  \item A model $\Phi \in X_5(K)$ is non-singular if and only if
    $\Delta(\Phi) \not= 0$.
  \item There exist $a_1,a_2,a_3,a_4,a_6 \in \Z[X_5]$ and $b_2,b_4,b_6
    \in \Z[X_5]$ satisfying
    \begin{equation}
      \label{tf}
      \begin{aligned}
        b_2 & = a_1^2+4 a_2, & b_4 & = a_1 a_3 + 2 a_4,
        \quad \quad b_6  =  a_3^2 + 4 a_6, \\
        c_4 & = b_2^2-24 b_4, \quad & c_6 & = -b_2^3 + 36 b_2 b_4 -
        216 b_6.
      \end{aligned}
    \end{equation}
  \item If $\Phi \in X_5(K)$ is non-singular then $\CC_\Phi$ has
    Jacobian elliptic curve
    \begin{equation*}
      y^2 + a_1 xy + a_3 y
      = x^3 + a_2 x^2 + a_4 x + a_6
    \end{equation*}
    where $a_i = a_i(\Phi)$.
  \end{enumerate}
\end{Theorem}
For the proof of Theorem~\ref{thm:invjac}(ii) we use the following lemma.
\begin{Lemma}
\label{lem:agen}
Let $c_4,c_6,\Delta \in R = \Z[x_1, \ldots,x_N]$ be primitive
polynomials satisfying $c_4^3 - c_6^2 = 1728 \Delta$. If there exists
$a_1 \in R$ satisfying $a_1^2 c_4 + c_6 \equiv 0 \pmod{4}$ then there
exist $a_2,a_3,a_4,a_6,b_2,b_4,b_6 \in R$ satisfying~(\ref{tf}).
\end{Lemma}

\begin{Proof}
  By unique factorisation in $\F_3[x_1, \ldots, x_N]$ and the 
  Chinese Remainder Theorem
  there exists $b_2 \in R$ with $c_4 \equiv b_2^2 \pmod{3}$,
  $c_6 \equiv -b_2^3 \pmod{3}$ and $b_2 \equiv a_1^2 \pmod{4}$.
  Then $b_2 c_4 + c_6 \equiv 0 \pmod{12}$ and $c_4^3 \equiv c_6^2 
  \equiv b_2^2 c_4^2 \pmod{24}$. Since $c_4$ is primitive it follows that
  $c_4 \equiv b_2^2 \pmod{24}$. Next putting $x = b_2$ in an 
  identity of Kraus \cite{Kraus},
  \[ (x^2 - c_4)^3 = ( x^3 - 3 x c_4 - 2 c_6) (x^3 + 2 c_6) + 3(x c_4 + c_6)^2
   + c_6^2 - c_4^3, \]
  we deduce $b_2^3 - 3 b_2 c_4 - 2 c_6 \equiv 0 \pmod{432}$. We put
  $b_4 = (b_2^2 - c_4)/24$ and $b_6 = (b_2^3 - 3 b_2 c_4 - 2 c_6)/432$.
  Then $0 \equiv c_4^3 - c_6^2 \equiv 16 b_2^2(b_2 b_6 - b_4^2) \pmod{64}$
  and so $b_2 b_6 \equiv b_4^2 \pmod{4}$. 
  By unique factorisation in $\F_2[x_1, \ldots, x_N]$ there exists
  $a_3 \in R$ with $b_4 \equiv a_1 a_3 \pmod{2}$. Then $b_4^2 
  \equiv a_1^2 a_3^2 \pmod{4}$ and $b_6 \equiv a_3^2 \pmod{4}$. 
  We put $a_2 = (b_2 - a_1^2)/4$, $a_4 = (b_4 - a_1 a_3)/2$ 
  and $a_6 = (b_6 - a_3^2)/4$. 
\end{Proof}

\begin{ProofOf}{Theorem~\ref{thm:invjac}}
  (i) This is \cite[Theorem 4.4(ii)]{g1inv}. \\
  (ii) By Lemma~\ref{lem:agen} it suffices to construct $a_1 \in
  \Z[X_5]$ with $a_1^2 c_4 + c_6 \equiv 0 \pmod{4}$.  In
  \cite[Section~10]{g1inv} we constructed an invariant $a_1 \in
  \F_2[X_5]$ of weight~1 and showed that together with $\Delta$ it
  generates the ring of invariants in characteristic~$2$.  In
  particular $c_4 \equiv a_1^4 \pmod{2}$ and $c_6 \equiv a_1^6
  \pmod{2}$.  So if we lift $a_1$ to $\Z[X_5]$ then $a_1^2 c_4 + c_6 =
  2 f$ for some $f \in \Z[X_5]$. Since $a_1$ is an invariant mod 2,
  $a_1^2$ is an invariant mod 4, and $f$ is an invariant mod 2.
  Therefore $f \equiv \lambda a_1^6 \pmod{2}$ for some $\lambda \in
  \{0,1\}$. Hence $a_1^2 c_4 \pm c_6 \equiv 0 \pmod{4}$.  Specialising
  to one of the Weierstrass models
  in~\cite[Section~6]{g1inv} shows that the sign is $+$. \\
  (iii) It is shown in \cite[Theorem 4.4(iii)]{g1inv} 
  that if $K$ is a perfect field with characteristic not $2$ or $3$ 
  then $\CC_{\Phi}$ has Jacobian $y^2 = x^3 -
  27 c_4(\Phi) x - 54 c_6(\Phi)$.  The proof is now identical to that
  of \cite[Theorem~2.10]{CFS}.  This generalises a result of Artin,
  Rodriguez-Villegas and Tate \cite{ARVT} in the case $n=3$.
\end{ProofOf}

\section{Minimisation theorems}
\label{sec:T}

Let $K$ be a discrete valuation field, with ring of integers $\OK$,
and normalised valuation $v : K^\times \to \Z$.  We assume throughout
that the residue field $k$ is perfect.  A genus one model $\Phi \in
X_5(K)$ is {\em integral} if it has coefficient in $\OK$. If $\Phi$ is
non-singular and integral then, by Theorem~\ref{thm:invjac} and the
standard formulae for transforming Weierstrass equations, we have
$v(\Delta(\Phi)) = v(\Delta_E) + 12 \level(\Phi)$ where $\Delta_E$ is
the minimal discriminant of $E = \Jac(\CC_{\Phi})$ and $\level(\Phi)$
is a non-negative integer we call the {\em level}.  We say that $\Phi$
is {\em minimal} if $v(\Delta(\Phi))$, or equivalently the level, is
minimal among all integral models $K$-equivalent to $\Phi$. Notice
that if $\Phi' = g \Phi$ for some $g=[A,B]$ with $A,B \in \GL_5(K)$
then $\level(\Phi') = \level(\Phi) + v(\det g)$.

\begin{Theorem}
\label{minthm}
Let $\Phi \in X_5(K)$ be non-singular.
\begin{enumerate}
\item (Weak minimisation theorem) If $\CC_{\Phi}(K) \not= \emptyset$
  then $\Phi$ is $K$-equivalent to an integral model of level $0$.
\item (Strong minimisation theorem) If $\CC_{\Phi}(L) \not= \emptyset$
  where $L$ is an unramified extension of $K$ then $\Phi$ is
  $K$-equivalent to an integral model of level $0$.
\end{enumerate}
\end{Theorem}

In this section we prove the weak minimisation theorem. 
In Section~\ref{sec:minalg} we describe an explicit algorithm 
for minimising. Inspection of this algorithm shows that
the minimal level is unchanged by an unramified extension.  
Theorem~\ref{minthm}(ii) then follows from Theorem~\ref{minthm}(i).
In Section~\ref{sec:insol} 
we prove a converse to the strong minimisation
theorem, thereby showing this result is best possible.

We refer to \cite[Section 2]{CFS} for notation and results analogous 
to those in
Section~\ref{g1mod} for genus one models of degree $4$, i.e.  quadric
intersections.  Let $E$ be an elliptic curve over $K$, and $D$
a $K$-rational divisor on $E$ of degree $n=4$ or $5$. The complete
linear system $|D|$ defines an embedding $E \subset \PP^{n-1}$.  The
image is defined by a genus one model~$\Phi\in X_n(K)$, and this model
is uniquely determined, up to $K$-equivalence, by the pair $(E,[D])$.
Moreover every non-singular model $\Phi \in X_n(K)$ with
$\CC_{\Phi}(K) \not= \emptyset$ arises in this way. Therefore
Theorem~\ref{minthm}(i) is an immediate consequence of the following.

\begin{Theorem}
\label{thmstar}
Let $E/K$ be an elliptic curve, with integral Weierstrass equation
\begin{equation}
\label{weqnstar}
y^2 + a_1 xy + a_3 y = x^3 + a_2 x^2 + a_4 x + a_6, 
\end{equation}
and let $D \in \Div_K(E)$ be a divisor on $E$ of degree $n=4$ or $5$.
Then $(E,[D])$ can be represented by an integral genus one model with
the same discriminant as~{\rm(\ref{weqnstar})}.
\end{Theorem}

The case $n=4$ is proved in \cite[Theorem 3.8]{CFS}.
To deduce the case $n=5$ from the case $n=4$ we use the following lemma.

\begin{Lemma}
\label{lem:4<->5}
Let $D \in \Div_K(E)$ be a divisor of degree $4$ and let 
$P \in E(K)$. Let $\ell_i, \alpha_i, \beta_i$ for $i=1,2,3$ 
be linear forms in $x_1, \ldots, x_4$ over $K$. 
The following statements are equivalent.
\begin{enumerate}
\item The pair $(E,[D])$ is represented by the quadric intersection
\begin{equation}
\label{eqnC4a}
 \begin{aligned}
   \ell_1 \alpha_1 + \ell_2 \alpha_2 + \ell_3 \alpha_3 &= 0 \\
   \ell_1 \beta_1 + \ell_2 \beta_2 + \ell_3 \beta_3 &= 0
 \end{aligned}
\end{equation}
and $P$ is the point defined by $\ell_1 = \ell_2 = \ell_3 = 0$.
\item The pair $(E,[D+P])$ is represented by the genus one model of degree $5$
\begin{equation}
\label{eqnC5}
\begin{pmatrix} 
  0 & \gamma & \alpha_1 & \alpha_2 & \alpha_3 \\
  & 0 & \beta_1 & \beta_2 & \beta_3 \\
  &   & 0 & \ell_3 & -\ell_2  \\
  &  - &   &  0 & \ell_1 \\
  & & & & 0
\end{pmatrix}
\end{equation}
where $\gamma = x_5$ and $P$ is the point $(x_1: \ldots :x_5) =
(0:\ldots:0:1)$.
\end{enumerate}
\end{Lemma}

\begin{Proof}
  An isomorphism $\psi: C_4 \to C_5$, between the curves $C_4$ and
  $C_5$ defined by~(\ref{eqnC4a}) and~(\ref{eqnC5}), is given by
  \[ \psi: (x_1:x_2:x_3:x_4) \mapsto (x_1 \ell_i : x_2 \ell_i : x_3
  \ell_i : x_4 \ell_i : \alpha_j \beta_k - \alpha_k \beta_j ) \]
  (where $i,j,k$ are any cyclic permutation of $1,2,3$) with inverse
  \[ \psi^{-1} : (x_1:x_2:x_3:x_4:x_5) \mapsto (x_1:x_2:x_3:x_4). \]
  This isomorphism identifies the points $\{\ell_1 =\ell_2 = \ell_3 =
  0\} \in C_4(K)$ and $(0:\ldots:0:1) \in C_5(K)$. To prove the
  equivalence of (i) and (ii) we note that if $C_4 \subset \PP^3$
  meets some plane in the divisor $D = P_1 + P_2 +P_3 + P_4$ then the
  points $\psi(P_i)$ and $(0:\ldots:0:1)$ are a hyperplane section for
  $C_5 \subset \PP^4$.
\end{Proof}

\begin{Lemma} 
\label{lem:same}
The genus one models~(\ref{eqnC4a}) and~(\ref{eqnC5})
have the same invariants. 
\end{Lemma}

\begin{Proof}
  Let $\Phi$ be the matrix~(\ref{eqnC5}) and write $P = \Pf(\Phi) =
  (p_1,\ldots, p_5)$. According to \cite[Section 5.4]{g1inv} there are
  invariant differentials $\omega_4$ on $C_4 = \{ p_1 = p_2 = 0\}
  \subset \PP^3$ and $\omega_5$ on $C_5 = \{ \rank \Phi \le 2 \}
  \subset \PP^4$ given by
 \[ \omega_n = \frac{x_1^2 d(x_2/x_1)}{\Omega_n(x_1,\ldots,x_n)} \]
where 
\[ \Omega_4 = \frac{\partial p_1}{\partial x_3} \frac{\partial
  p_2}{\partial x_4} - \frac{\partial p_1}{\partial x_4}
\frac{\partial p_2}{\partial x_3} \qquad \text{ and } \qquad \Omega_5
= \frac{\partial P}{\partial x_3} \,\, \frac{\partial \Phi}{\partial
  x_5} \,\, \frac{\partial P^T}{\partial x_4}.  \] 
In the expression for $\Omega_5$ we have written the partial
derivative of a matrix as a short-hand for the matrix of
partial derivatives.  Since the only entries of $\Phi$ to involve
$x_5$ are in the top left $2 \times 2$ submatrix, it is clear that $\Omega_4
= \pm \Omega_5$. Hence the isomorphism $\psi : C_4 \to C_5$ identifies
the invariant differentials $\omega_4$ and $\omega_5$ (up to sign). It
follows by~\cite[Proposition 5.23]{g1inv} that~(\ref{eqnC4a})
and~(\ref{eqnC5}) have the same invariants $c_4$, $c_6$ and $\Delta$.
\end{Proof}

\begin{ProofOf}{Theorem~\ref{thmstar}}
Let $D \in \Div_K(E)$ be a divisor of degree 4, and let
$P \in E(K)$. We show that if the theorem holds for $D$
then it holds for $D+P$. Suppose $(E,[D])$ is represented by 
an integral quadric intersection
with discriminant $\Delta$. Since $\OK$ is a principal ideal domain, 
$\SL_4(\OK)$ acts transitively on $\PP^3(K)$.
So we may assume $P$ is the point $(x_1: x_2 :x_3:x_4 )= (0:0:0:1)$. 
Our model is now of the form~(\ref{eqnC4a}) with $\ell_i =x_i$
for $i=1,2,3$. We may choose the linear forms $\alpha_i$ and $\beta_i$
to have coefficients in $\OK$. Then the genus one model~(\ref{eqnC5}) 
is an integral model of discriminant $\Delta$ representing 
the pair $(E,[D+P])$. 
\end{ProofOf}

\section{Minimisation algorithms}
\label{sec:minalg}
For $\Phi \in X_5(\OK)$ we write $\phi \in X_5(k)$ for its reduction 
mod~$\pi$.  The {\em singular~locus}
$\Sing \CC_\phi$ is the set of points $P \in \CC_\phi$ with tangent space
of dimension greater than~$1$.  (We make this definition regardless of
whether $\CC_\phi$ is a curve. In particular 
all points on components of dimension
at least $2$ are singular.) For example, if $\phi$ takes the
form~(\ref{eqnC5}) with $\gamma = x_5$ and $\ell_i, \alpha_i,
\beta_i$ linear forms in $x_1, \ldots, x_4$, then
$P=(0:\ldots:0:1)$ is singular if and only if $\ell_1, \ell_2, \ell_3$
are linearly dependent.  An integral genus one model $\Phi \in
X_5(\OK)$ is {\em saturated} if its $4 \times 4$ Pfaffians $p_1, \ldots,p_5$
are linearly independent mod $\pi$. We write $I_m$ for the $m \times m$
identity matrix.

Our algorithm for minimising genus one models of degree $5$
generalises the algorithm for models of degree $3$ in
\cite[Section 4B]{CFS}.

\begin{Theorem} 
\label{minalgthm}
Let $\Phi \in X_5(\OK)$ be saturated and of positive level. 
\begin{enumerate}
\item The singular locus $\Sing \CC_\phi$ 
  does not span $\PP^{4}$.
\item Let $B \in \GL_5(\OK)$ represent a change of co-ordinates on
  $\PP^4$ mapping the linear span of the singular locus in (i) to 
  $\{ x_{m+1} = \ldots = x_{5} = 0 \}$.  Then there exist $A \in
  \GL_5(K)$ and $\mu \in K^\times$ such that 
  $[ A , \mu \Diag(I_m, \pi I_{5-m})B] \Phi$ is an integral model of
  the same or smaller level.
\item If $\Phi$ is non-minimal then repeating the procedure in (ii)
  either gives a non-saturated model or decreases the level after
  finitely many iterations.
\end{enumerate}
\end{Theorem}

As it stands Theorem~\ref{minalgthm} does not give an algorithm for
minimising since we must show how to find $A$ and $\mu$ in (ii), and
show how to decrease the level of a non-saturated model.  We do this
in Theorem~\ref{thm:sat} below. Theorem~\ref{minalgthm} 
is proved in Sections~\ref{sec:sing} and~\ref{sec:wtsl}. 
In Section~\ref{sec:numit} 
we bound the number of iterations required in~(iii).

\begin{Theorem}
\label{thm:sat}
Let $\Phi \in X_5(\OK)$ be non-singular. Let $\ell_0$ be the minimum
of the levels of all integral models that are $K$-equivalent to $\Phi$
via a transformation of the form $[A, \mu I_5]$ where $A \in \GL_5(K)$
and $\mu \in K^\times$. 
\begin{enumerate} \item 
We may compute an integral model of the form $[A, \mu I_5] \Phi$ 
with level~$\ell_0$ as follows:
\begin{description}
\item[Step 1] Write $\Pf(\Phi) = (p_1, \ldots, p_5)$. Compute 
  $A =(a_{ij})\in \GL_5(K)$ and quadrics 
  $q_1, \ldots, q_5 \in \OK[x_1, \ldots,x_5]$
  such that $p_j = \sum_{i=1}^5 a_{ij} q_i$ and $q_1, \ldots,q_5$ are
  linearly independent mod $\pi$.  Then replace $\Phi$ by $[A,\mu I_5]
  \Phi$ where $\mu \in K^\times$ is chosen so that $\Phi$ has
  coefficients in $\OK$ not all in $\pi \OK$.
\item [Step 2] Replace $\Phi$ by  $[A,I_5] \Phi$ where $A \in \GL_5(\OK)$
  is chosen so that the first two rows of $\Phi$ are divisible by 
  $\pi^e$, with $e \ge 0$ as large as possible. 
  Then divide the first row and column by $\pi^e$.
\end{description}
\item 
If the model computed in Step~$1$ is non-saturated, then 
we may compute an integral model of level smaller 
  than $\ell_0$ by modifying Step $2$ so that we divide the 
  first two rows and columns 
  by $\pi^e$, and then make a transformation of the form $[I_5,B]$ 
  to preserve integrality.
\end{enumerate}
\end{Theorem}

\begin{Proof}
  With the notation of Step 1 we have
  \[ \Pf(A \Phi A^T) = \Pf(\Phi) \adj A = (q_1, \ldots, q_5) A \adj A
  = (\det A) (q_1, \ldots, q_5). \] 
  So after Step 1 we have $\Pf(\Phi) = (\la q_1, \ldots, \la q_5)$
  where $\la := \mu^2 \det A \in \OK$. We split into the cases
  $v(\la)=0$ and $v(\la) \ge 1$. First we need two lemmas.

  \begin{Lemma}
    \label{lem:sat}
    Let $\Phi, \Phi' \in X_5(\OK)$ be non-singular models with
    $\Phi' = [A, \mu I_5] \Phi$ for some $A \in \GL_5(K)$ and $\mu \in
    K^\times$.
    \begin{enumerate}
    \item If $\Phi$ is saturated then $\level(\Phi') \ge \level(\Phi)$
      with equality if and only if $\Phi$ and $\Phi'$ are
      $\OK$-equivalent.
    \item If $\Phi$ and $\Phi'$ are of the form output by Step $1$ 
      then they are $\OK$-equivalent.
    \end{enumerate}
  \end{Lemma}
  \begin{Proof}
    We have $\Pf(\Phi') = \Pf(\Phi) M$ where $M := \mu^2 \adj A$. \\
    (i) Since $\Phi$ is saturated, $M$ has entries in $\OK$.  Hence
    $\level(\Phi') - \level(\Phi) = \tfrac{1}{2} v(\det M) \ge 0$ with
    equality if and only if $M \in \GL_5(\OK)$. If $M \in \GL_5(\OK)$
    then replacing $[A,\mu I_5]$ by $[\lambda A,\lambda^{-2} \mu I_5]$
    for suitable $\lambda \in K^\times$ we may assume $A \in
    \GL_5(\OK)$.  Since $\Phi$ and $\Phi'$ have the same level they must
    therefore be $\OK$-equivalent. \\
    (ii) Since $\Pf(\Phi)$ and $\Pf(\Phi')$ are scalar multiples of
    bases for the same $\OK$-module, some scalar multiple of $M$
    belongs to $\GL_5(\OK)$. Replacing $[A,\mu I_5]$ by $[\lambda
    A,\lambda^{-2} \mu I_5]$ for suitable $\lambda \in K^\times$ we
    may assume $A \in \GL_5(\OK)$.  Since $\Phi$ and $\Phi'$ are
    primitive they must therefore be $\OK$-equivalent. 
  \end{Proof}

  \begin{Lemma}
   \label{lem:pfzero}
  Let $\phi \in X_5(k)$ be a genus one model all of whose $4 \times 4$
  Pfaffians are identically zero. Then $\phi$ is $k$-equivalent to either
  \[
  \begin{pmatrix}
    0 & \ell_2 & \ell_3 & \ell_4 & \ell_5 \\
    & 0 & 0 & 0 & 0 \\
    & & 0 & 0 & 0 \\
    & - & & 0 & 0 \\
    & & & & 0
  \end{pmatrix}
  \quad \text { or } \quad
  \begin{pmatrix}
    0 & x_1 & x_2 & 0 & 0 \\
    & 0 & x_3 & 0 & 0 \\
    & & 0 & 0 & 0 \\
    & - & & 0 & 0 \\
    & & & & 0
  \end{pmatrix} \] where $\ell_2, \ldots, \ell_5$ are linear forms.
\end{Lemma}

\begin{Proof} This is clear.
\end{Proof}

We now complete the proof of Theorem~\ref{thm:sat}. If $v(\la) = 0$ then
$\Phi$ is saturated and we are done by Lemma~\ref{lem:sat}(i).  
So suppose $e := v(\la) \ge 1$. In Step~1 the matrix $A$ has entries 
in $\OK$. So $v(\mu) \le 0$ and the level is increased by 
\[2 v(\det A) + 5 v(\mu) \le 2 v(\mu^2 \det A) = 2 e.\]  
Lemma~\ref{lem:sat}(ii) shows that when we apply
Step~1 to both $\Phi$ and the model implicit in the definition of
$\ell_0$ then we obtain models that are $\OK$-equivalent. So it
will suffice to show that Step~2 reduces the level by $2e$, whereas the
modified version in (ii) reduces the level by more than $2e$.

Since $\Pf(\Phi) = (\la q_1, \ldots, \la q_5)$ we have $(q_1, \ldots,
q_5) \Phi = 0$.  The reduction of $\Phi$ takes one of the forms
specified in Lemma~\ref{lem:pfzero}.  In the first case we have $q_1
\ell_{j} \equiv 0 \pmod{\pi}$ for $j=2,\ldots,5$. This contradicts the
choices of $q_1, \ldots, q_5$ and $\mu$ in Step~1.  So we must be in the 
second case. 
Replacing $\Phi$ by an $\OK$-equivalent model we may assume it
takes the form~(\ref{eqnC5}) with $\ell_i = x_i$ for $i=1,2,3$ and
$\alpha_1, \alpha_2, \alpha_3, \beta_1, \beta_2, \beta_3, \gamma$ 
linear forms that vanish mod $\pi$.  By row and column operations we
may assume $\alpha_2 \in \langle x_2 , \ldots, x_5 \rangle$ and
$\alpha_3 \in \langle x_3, \ldots, x_5 \rangle$.  Then since $\pi^e
\dv (x_1 \alpha_1 + x_2 \alpha_2 + x_3 \alpha_3)$ we have $\pi^e \dv
\alpha_1, \alpha_2, \alpha_3$.  Likewise we may assume $\pi^e \dv
\beta_1, \beta_2, \beta_3$.  The remaining Pfaffians show that $\pi^e
\mid \gamma$. Steps $2$ and its modified version in (ii) now reduce
the level by $2e$ and $3e$ respectively.
\end{Proof}

\begin{Corollary}
\label{cor:wlog}
For the proof of Theorem~\ref{minalgthm} we are free to replace $\Phi$
by an $\OK$-equivalent model, and to replace $K$ by an unramified
field extension.
\end{Corollary}
\begin{Proof} Let $\Phi_1, \Phi_2 \in X_5(\OK)$ be $\OK$-equivalent
  models and $\Phi'_1, \Phi'_2 \in X_5(\OK)$ the models returned by
  Theorem~\ref{minalgthm}(ii). Lemma~\ref{lem:sat}(i) and 
  \cite[Lemma 4.1]{CFS} together show that if $\Phi'_1$ is saturated and
  $\level(\Phi'_1) = \level(\Phi'_2)$ then $\Phi'_1$ and $\Phi'_2$ are
  $\OK$-equivalent.  Therefore the number of iterations required in
  Theorem~\ref{minalgthm}(iii) depends only on the $\OK$-equivalence
  class of $\Phi$.

  For the final statement we note that the performance of the algorithms 
  in Theorems~\ref{minalgthm} and~\ref{thm:sat} is unchanged by an
  unramified field extension.
\end{Proof}

Replacing $K$ by its strict Henselisation, we may assume in the next 
three sections that $K$ is Henselian and its residue field $k$ 
is algebraically closed.

\section{The singular locus}
\label{sec:sing}

In this section and the next we prove Theorem~\ref{minalgthm}.
\begin{Lemma} 
\label{lem:null}
Let $\phi \in X_5(k)$ be a genus one model. Suppose $\Gamma
\subset \CC_\phi$ is either a line or a (non-singular) conic. Then either
$\Gamma \subset \Sing \CC_\phi$ or
\[  \# (\Gamma \cap \Sing \CC_\phi)  = \left\{ \begin{array}{ll} 
1 & \text{ if } c_4(\phi) = c_6(\phi) = 0, \\
2 & \text{ otherwise. } \end{array} \right. \]
\end{Lemma}

\begin{Proof}
  (i) If $\CC_\phi$ contains the line $\Gamma = \{x_3 = x_4 = x_5 = 0
  \}$, but not every point on $\Gamma$ is singular, then (unless
  $\CC_\phi$ is a cone -- which is an easy special case with 
$c_4(\phi) = c_6(\phi) = 0$) we may suppose
  $\phi$ is
  \begin{equation*}
    \begin{pmatrix}
      0  & x_1 & x_2 & * & * \\
      &  0  & * &\alpha &\beta\\
      &     &  0  &\gamma&  \delta  \\
      &  -  &     &  0  &  x_5  \\
      & & & & 0
    \end{pmatrix}
  \end{equation*}
  where $\alpha,\beta,\gamma,\delta$ 
  and the entries $*$ are linear forms in $x_3,x_4,x_5$. By row and column
  operations (and substitutions for $x_1$ and $x_2$) we may suppose
  $\alpha,\beta,\gamma,\delta$ do not involve $x_5$.
 We write $\alpha
  = \alpha_3 x_3 + \alpha_4 x_4, \ldots, \delta = \delta_3 x_3 +
  \delta_4 x_4$ and put
  \[ q(s,t) = \det \left(
    \begin{pmatrix} \gamma_3 & \gamma_4 \\ \delta_3 & \delta_4
    \end{pmatrix} s - \begin{pmatrix} \alpha_3 & \alpha_4 \\ \beta_3
      & \beta_4
    \end{pmatrix} t \right). \] By the Jacobian criterion we have
  \[ \Gamma \cap \Sing \CC_\phi = \{ (s:t:0:0:0) \mid q(s,t) = 0 \}. \]
  A calculation using Lemma~\ref{lem:same} shows that $c_4(\phi) =
  \Delta(q)^2$ and $c_6(\phi) = -\Delta(q)^3$ where $\Delta(q)$ is the
  discriminant of the binary quadratic form $q$.

  (ii) Suppose $\CC_\phi$ contains the conic $\Gamma = \{
  f(x_1,x_2,x_3) = x_4 = x_5 = 0 \}$, but not every point on $\Gamma$
  is singular. Let $\Pf(\phi) = (p_1, \ldots,p_5)$.  Replacing $\phi$
  by an equivalent model we may suppose $p_i(x_1,x_2,x_3,0,0)=0$ for
  $i=1,2,3,4$ and $p_5(x_1,x_2,x_3,0,0)=f$.  Since
  $\Pf(\phi) \phi = 0$, and $\Gamma$ is not contained in any component
  of $\CC_\phi$ of higher dimension, 
  we may further suppose the last column of $\phi$
  has entries $x_4,x_5,0,0,0$.  The monomials appearing in the invariants
  $c_4$ and $c_6$ are limited by the fact they are invariant under all 
  pairs of diagonal matrices.  These restrictions show that $c_4(\phi)$ and
  $c_6(\phi)$ are unchanged if we set $x_4=x_5=0$ in all entries of
  $\phi$ outside the last row and column.  Writing $f = \sum_{i \le j} 
  a_{ij} x_i x_j$ and $\phi_{34} = \sum b_i x_i$ we put
  \[\delta = \left| \begin{array}{cccc}
      2 a_{11} & a_{12} & a_{13} & b_1 \\ a_{12} & 2 a_{22} & a_{23} & b_2 \\
      a_{13} & a_{23} & 2 a_{33} & b_3 \\ b_1 & b_2 & b_3 & 0 
    \end{array} \right|. \]
  A calculation using Lemma~\ref{lem:same} shows that $c_4(\phi) = \delta^2$
  and $c_6(\phi) = -\delta^3$. By a change of co-ordinates we may
  suppose $f = x_1 x_3 - x_2^2$. Then $\delta$ is the discriminant of the
  binary quadratic form $q(s,t) = \phi_{34}(s^2,st,t^2,0,0)$ and  
  by the Jacobian criterion 
  \[ \Gamma \cap \Sing \CC_\phi = \{ (s^2:st:t^2:0:0) \mid q(s,t) = 0 \}.
  \vspace{-3ex} \]
\end{Proof}

\begin{Lemma}
\label{lem:sing}
Let $\phi \in X_5(k)$ be a genus one model.  Suppose the $4 \times 4$
Pfaffians $p_1, \ldots, p_5$ are linearly independent and $c_4(\phi) =
c_6(\phi) = 0$. Then either $\Sing \CC_\phi$ is a linear subspace of
$\PP^4$ or $\phi$ is equivalent to a model of the form
\begin{equation}
\label{x45}
 \begin{pmatrix}
  0  & \xi &  \alpha & \beta & \eta \\
     &  0  &\gamma &\delta &x_5\\
     &     &  0  &x_5&  0  \\
     &  -  &     &  0  &  0  \\
     &     &     &     &  0  
\end{pmatrix} 
\end{equation}
where $\xi,\eta, \alpha, \beta, \gamma, \delta$ are linear forms in
$x_1, \ldots, x_5$.
\end{Lemma}

\begin{Proof} If $P_1, P_2 \in \Sing \CC_\phi$ are distinct and the line
  $\ell$ between them is contained in $\CC_\phi$ then by
  Lemma~\ref{lem:null} we have $\ell \subset \Sing \CC_\phi$.  So either
  $\Sing \CC_\phi$ is a linear subspace of $\PP^4$ or there exist $P_1,
  P_2 \in \Sing \CC_\phi$ joined by a line not contained in $\CC_\phi$. We
  move these points to $(1:0:\ldots:0)$ and $(0:1:\ldots:0)$. 
  Writing $\phi = \sum x_i M_i$, the matrices $M_1$ and $M_2$ have
  rank $2$ but their sum has rank $4$. Therefore $\phi$ is equivalent
  to a model with $\phi_{12} = x_1$, $\phi_{34} = x_2$ and all other
  $\phi_{ij}$ (for $i<j$) linear forms in $x_3,x_4,x_5$.  Since $P_1$
  and $P_2$ are singular, $\phi_{35}$ and $\phi_{45}$ are linearly
  dependent, and $\phi_{15}$ and $\phi_{25}$ are linearly
  dependent. So the space of linear forms spanned by the entries of
  the last column has dimension at most $2$. In fact it has dimension
  exactly $2$, since $p_1, \ldots,p_5$ are linearly independent.

  Replacing $\phi$ by an equivalent model we may assume it has last
  column with entries $x_4,x_5,0,0,0$. The transformation used here
  does not move $P_1$ and $P_2$ but may change the matrices $M_1$ and
  $M_2$.  Let $\Gamma = \{ x_4 = x_5 = p_5 = 0 \} \subset \CC_\phi$.
  Then $P_1$ and $P_2$ are contained in $\Gamma$ but the line between
  them is not. It follows that $\Gamma$ is either a non-singular conic
  or a pair of concurrent lines. In either case Lemma~\ref{lem:null}
  shows that $\Gamma \subset \Sing \CC_\phi$.  By the Jacobian criterion
  it follows that $\phi_{34} \in \langle x_4, x_5 \rangle$. However
  $\phi_{34}$ is non-zero since $p_1, \ldots,p_5$ are linearly
  independent. Therefore $\phi$ is equivalent to a model of the
  form~(\ref{x45}).
\end{Proof}

\begin{Lemma}
\label{lem:w34}
Let $\Phi \in X_5(\OK)$ be a saturated non-singular model with
reduction $\phi$ of the form~(\ref{x45}). Suppose $\Sing \CC_\phi$ has
linear span $\{ x_{m+1} = \ldots = x_5= 0\}$.
\begin{enumerate}
\item There exist $A \in \GL_5(K)$ and $\mu \in K^\times$ such that $[
  A , \mu \Diag(I_m, \pi I_{5-m})] \Phi$ is an integral model of the
  same or smaller level.
\item Suppose that either $\delta= 0$ and $\Phi_{45} \equiv 0 
\pmod{\pi^2}$, or $\Phi_{35} \equiv \Phi_{45} \equiv 0  \pmod{\pi^2}$.
Then there is a transformation as in (i) that decreases the level.
\end{enumerate}
\end{Lemma}
\begin{Proof}
Computing the $4 \times 4$ Pfaffians of~(\ref{x45}) we find
\begin{equation}
\label{conicandline}
\CC_\phi =  \{ \eta = x_5 = \alpha \delta - \beta \gamma =0 \} 
\cup \{ \gamma = \delta = x_5 = 0 \}. 
\end{equation}

First suppose $\gamma, \delta, x_5$ are linearly dependent. By an
$\OK$-equivalence we may assume $\delta = 0$. Then $\{ \gamma = x_5 =
0 \} \subset \Sing \CC_{\phi} \subset \{ x_5 = 0 \}$.  Therefore $m=3$
or $4$. The required transformations are as follows.
\[ \begin{array}{r|c|c} & m = 3 & m = 4 \\ \hline {\text{(i)}} & A=
  \Diag(\pi,1,1,1,1), \quad \mu = \pi^{-1} &
  A= \Diag(\pi,\pi,1,1,1), \quad \mu = \pi^{-1} \\
  {\text{(ii)}} & A= \Diag(\pi,1,1,1,1), \quad \mu = \pi^{-1} & A =
  \Diag(\pi,1,1,\pi^{-1},\pi^{-1}), \quad \mu = 1
\end{array} \]

Now suppose $\gamma, \delta, x_5$ are linearly independent. Since
$\Phi$ is saturated $\eta, x_5$ are linearly independent. A
calculation shows that $\Sing \CC_\phi$ is the first of the two
components in~(\ref{conicandline}). Therefore $m=2$ or $3$.  If $m=2$
then we may assume $\beta, \gamma, \delta, \eta$ are linear
combinations of $x_3, x_4, x_5$. The required transformations are as
follows
\[ \begin{array}{r|c|c} & m = 2 & m = 3 \\ \hline {\text{(i)}} & A=
  \Diag(\pi,1,1,1,1), \quad \mu = \pi^{-1} &
  A= \Diag(1,1,1,1,\pi^{-1}), \quad \mu = 1 \\
  {\text{(ii)}} & A = \Diag(1,1,1,\pi^{-1},\pi^{-1}), \quad \mu = 1 &
  A = \Diag(\pi,\pi,1,1,\pi^{-1}), \quad \mu = \pi^{-1}
\end{array} \]
\end{Proof}

We now prove the first two parts of Theorem~\ref{minalgthm}.  Let
$\Phi \in X_5(\OK)$ be saturated and of positive level.
Lemma~\ref{lem:sing} shows that either $\Sing \CC_\phi$ is a linear
subspace 
or $\CC_\phi$ is contained in a hyperplane. 
Since $\CC_\phi$ is defined by $5$ linearly independent quadrics it 
cannot be all of $\PP^4$. This proves Theorem~\ref{minalgthm}(i).

The proof of Theorem~\ref{minalgthm}(ii) in the case $\phi$ takes the
form~(\ref{x45}) was already given in Lemma~\ref{lem:w34}(i).  So by
Lemma~\ref{lem:sing} we may assume $\Sing \CC_\phi = \{ x_{m+1} = \ldots
= x_5 = 0 \}$.  We apply Lemma~\ref{lem:pfzero} to the
reduction mod $\pi$ of $[I_5,\Diag(I_m,\pi I_{5-m})] \Phi$. In the
second case of that lemma we have $m \ge 3$. 
We take $A= \Diag(1,1,1,1,\pi^{-1})$ and $\mu = 1$.  
Otherwise we are in the first case. If $m \ge 2$ then we take $A=
\Diag(\pi,1,1,1,1)$ and $\mu = \pi^{-1}$.  It remains to treat the
case $m=1$, in other words the case $\Sing \CC_\phi$ is a point.

By \cite[Lemma 5.8]{g1inv} every component of
$\CC_\phi$ has dimension at least $1$. So if $\Sing \CC_\phi$ is
just a point then there are also smooth points on $\CC_\phi$.  Since $K$ is
Henselian it follows that $\CC_\Phi(K) \not= \emptyset$ and so, by
Theorem~\ref{minthm}(i), $\Phi$ is non-minimal. With this extra
hypothesis we show in the next section that the singular point on
$\CC_\phi$ is non-regular (as a point on the $\OK$-scheme $\CC_{\Phi}$).

We may suppose $\phi_{12} = x_1$ and all other $\phi_{ij}$ (for $i<j$)
are linear forms in $x_2, \ldots, x_5$. Since $P = (1:0: \ldots: 0)$
is singular, $\phi_{34}, \phi_{35}, \phi_{45}$ are linearly dependent.
So replacing $\Phi$ by an $\OK$-equivalent model we may assume
$\phi_{45}=0$.  In the presence of the stronger condition that $P$ is
non-regular we may further arrange that the coefficient of $x_1$ in
$\Phi_{45}$ is divisible by $\pi^2$. Taking $A=
\Diag(1,1,1,\pi^{-1},\pi^{-1})$ and $\mu = 1$ now preserves the level.

\section{Weights and slopes}
\label{sec:wtsl}

In this section we complete the proof of Theorem~\ref{minalgthm}.

\begin{Definition} (i) The set of {\em weights} is 
  \[ \W = \left\{ (r,s) \in \Z^{5} \times \Z^{5} {\Bigg|}
    \begin{array}{c} r_1 \le r_2 \le \ldots \le r_5, \quad
      s_1 \le s_2 \le \ldots \le s_5, \\
      2 \sum_{i=1}^5 r_i = 1 + \sum_{i=1}^5 s_i \end{array} \right\}.
  \] (ii) A {\em weight for} $\Phi \in X_5(\OK)$ is $(r,s) \in \W$ 
  such that the model
  \begin{equation}
    \label{eqn:trans}
    [\Diag(\pi^{-r_1}, \ldots, \pi^{-r_5}),\Diag(\pi^{s_1}, \ldots,
    \pi^{s_5}) ] \Phi 
\end{equation}
has coefficients in $\OK$. \\
(iii) Let $w = (r,s)$ and $w' = (r',s')$ be weights. Then $w$
{\em dominates} $w'$ if
\[ \max(r_i + r_j - s_k, 0) \ge \max(r'_i + r'_j - s'_k, 0) \] for all
$1 \le i < j \le 5$ and $1 \le k \le 5$.
\end{Definition}

Let ${\bf 1} = (1,1,\ldots, 1)$. Then $\la \in \Z$ acts on $\W$ as
$(r,s) \mapsto (r + \la {\bf 1},s + 2 \la {\bf 1})$. Since weights 
in the same $\Z$-orbit determine the same transformation~(\ref{eqn:trans}) 
we may regard such weights as equivalent.

\begin{Lemma}
\label{lem:wts}
Let $\Phi \in X_5(\OK)$ be an integral genus one model. \\
(i) If $\Phi$ is non-minimal then it is $\OK$-equivalent
to a model with a weight. \\
(ii) If $\Phi$ has weight $w$ and $w$ dominates $w'$ then $\Phi$ has
weight $w'$.
\end{Lemma}
\begin{Proof} (i) By hypothesis there exist $A,B \in \GL_5(K)$ with
  $[A,B] \Phi$ integral and $2 v(\det A) + v(\det B) =-1$.
  We put $A$ and $B$ in Smith normal form. \\
  (ii) Let $\Phi = (\Phi_{ij})$ with $\Phi_{ij} = \sum_{k} a_{ijk}
  x_k$.  Then $\Phi$ has weight $(r,s)$ if and only if $v(a_{ijk}) \ge
  \max( r_i + r_j - s_k,0)$ for all $1 \le i < j \le 5$ and $1 \le k
  \le 5$.
\end{Proof}

\begin{Lemma}
\label{lem:showsing}
Let $\Phi \in X_5(\OK)$ have weight $(r,s) \in \W$ with either $r_1+r_4>s_1$
or $r_2+r_3>s_1$. Then $P = (1:0:\ldots:0) \in \CC_\phi$ is a singular
point. Moreover if $s_1 < s_3$ then $P$ is non-regular (as a point on 
the $\OK$-scheme $\CC_{\Phi}$). 
\end{Lemma}
\begin{Proof}
  We write $\phi = \sum x_i M_i$. If $r_1 + r_4>s_1$ then the only
  non-zero entries of $M_1$ are in the top left $3 \times 3$
  submatrix. If $r_2 + r_3>s_1$ then the only non-zero entries of
  $M_1$ are in the first row and column. In both cases $\rank M_1
  \le 2$ and so $P \in \CC_\phi$.  If $M_1=0$ then $P$ is singular
  (and non-regular). So we may assume $M_1 \not= 0$. We are free to
  multiply rows of $\Phi$ by units in $\OK$ and to subtract
  $\OK$-multiples of later rows from earlier rows (it being understood
  that we also make 
  the corresponding column operations).  In particular these operations
  do not upset our hypothesis that $\Phi$ has weight $(r,s)$. Let
  $E_{ij}$ be the $5 \times 5$ matrix with a $1$ in the $(i,j)$-place
  and zeros elsewhere. By row and column operations we reduce to 
  the case $M_1 = E_{ij} - E_{ji}$ where $(i,j) \in \{
  (1,2),(1,3),(1,4), (1,5),(2,3)\}$. Let $a < b < c$ be chosen such
  that $\{i,j,a,b,c\} = \{1, \ldots, 5 \}$. Since $r_i + r_j\le s_1 \le
  s_2$ it follows by the definition of $\W$ that \[s_3 + s_4 + s_5 <
  (r_a + r_b) + (r_a + r_c) + (r_b + r_c).\] Therefore at least one of
  the following three inequalities holds:
  \[ \begin{aligned} s_3 & < r_a + r_b && \implies & \phi_{ab},
    \phi_{ac}, \phi_{bc} & \in
    \langle x_4, x_5 \rangle, \\
    s_4 & < r_a + r_c && \implies & \phi_{ac}, \phi_{bc} & \in
    \langle x_5 \rangle, \\
    s_5 & < r_b + r_c && \implies & \phi_{bc} & = 0.
  \end{aligned}
\]
Since the tangent space at $P$ is $\{ \phi_{ab}= \phi_{ac}= \phi_{bc}
= 0 \}$ it follows that $P \in \CC_\phi$ is a singular point.

If $s_1 < s_3$ then the same argument shows there is some
$\OK$-linear combination of $\Phi_{ab}, \Phi_{ac}, \Phi_{bc}$ (with
not all coefficients in $\pi \OK$) that not only vanishes mod $\pi$,
but whose coefficient of $x_1$ vanishes mod $\pi^2$. Hence 
$P$ is non-regular.
\end{Proof}

\begin{Lemma}
  \label{lem:sevenwts}
  Let $(r,s) \in \W$ be a weight with $r_1+r_4 \le s_1$ and $r_2+r_3
  \le s_1$.  Then $(r,s)$ dominates one of the weights $w_1, \ldots,
  w_7$ in the following table.  \small
\[ \begin{array}{c|ccccc|ccccc}
 & r_1 & r_2 & r_3 & r_4 & r_5 & s_1 & s_2 & s_3 & s_4 & s_5 \\ \hline
 w_1 & 0 & 0 & 0 & 0 & 1 & 0 & 0 & 0 & 0 & 1 \\
 w_2 & 0 & 0 & 1 & 1 & 1 & 1 & 1 & 1 & 1 & 1 \\
 w_3 & 0 & 0 & 1 & 1 & 2 & 1 & 1 & 1 & 2 & 2 \\
 w_4 & 0 & 1 & 1 & 2 & 2 & 2 & 2 & 2 & 2 & 3 \\
 w_5 & 0 & 1 & 1 & 2 & 3 & 2 & 2 & 2 & 3 & 4 \\
 w_6 & 0 & 1 & 1 & 2 & 3 & 2 & 2 & 3 & 3 & 3 \\
 w_7 & 0 & 1 & 2 & 3 & 4 & 3 & 3 & 4 & 4 & 5 
\end{array} \]
\normalsize
\end{Lemma}

\begin{Proof}
  We checked the lemma by writing a computer program using the
  simplex algorithm. 
 See the proof of Lemma~\ref{lem:wt29} for details.
\end{Proof}

\begin{Definition}
  The {\em slope} of $\Phi \in X_5(\OK)$ is the least possible value
  of $v(\det B)$ for $B \in \GL_5(K)$ a matrix with entries in $\OK$
  for which there exist $A \in \GL_5(K)$ and $\mu \in K^\times$ such
  that $[A, \mu B] \Phi$ is an integral model of smaller level.
\end{Definition}

We now complete the proof of Theorem~\ref{minalgthm}.  Since $\Phi \in
X_5(\OK)$ is non-minimal it has a slope $\sigma$, say. 
Lemma~\ref{lem:sat}(i) shows that if $\sigma = 0$ then $\Phi$ 
is non-saturated.  So we may
assume $\sigma > 0$. By Lemma~\ref{lem:wts} (and Corollary~\ref{cor:wlog})
we may replace $\Phi$ by
an $\OK$-equivalent model with a weight, say $(r,s)$. Moreover we may
assume the weight realises the slope, i.e.  $\sigma = \sum_{i=1}^5
(s_i - s_1)$.

Suppose that either $r_1+r_4>s_1$ or $r_2+r_3>s_1$.  Since $\sigma >
0$ there exists $1 \le m \le 4$ such that $s_1 = \ldots = s_m <
s_{m+1}$.  Lemma~\ref{lem:showsing} shows (by first making
unimodular transformations involving only $x_1, \ldots, x_m$) that
\begin{equation}
\label{linsing}
 \{ x_{m+1} = \ldots = x_{5} = 0 \} \subset \Sing \CC_\phi \enskip . 
\end{equation}
Moreover if $m=1$ then the point we have constructed is non-regular.
(This is needed to complete the proof of Theorem~\ref{minalgthm}(ii)
at the end Section~\ref{sec:sing}.)

Regardless of whether we have equality in~(\ref{linsing}) it follows
that if the level is preserved then the slope is decreased. So after
finitely many iterations $\Phi$ is either non-saturated or has weight
$(r,s)$ with $r_1+r_4 \le s_1$ and $r_2+r_3 \le s_1$. 
In this last
case Lemmas~\ref{lem:wts} and~\ref{lem:sevenwts} show that $\Phi$ has
weight $w$ for some $w \in \{ w_1, \ldots, w_7 \}$. If $w \in \{w_1,
w_2, w_6\}$ then $\Phi$ is non-saturated.  If $w \in \{w_5, w_7\}$
then $\Phi$ is $\OK$-equivalent to a model with weight $w_3$. 
(This is achieved by a unimodular transformation involving only the
second and third rows and columns, respectively a unimodular transformation
involving only $x_3$ and $x_4$.) Finally
if $w \in \{w_3, w_4\}$ then $\Phi$ is $\OK$-equivalent to a model 
of the form considered in Lemma~\ref{lem:w34}(ii) 

\section{The number of iterations}
\label{sec:numit}

We have shown that if we start with a non-minimal model
then iterating the procedure in Theorem~\ref{minalgthm}(ii) 
eventually gives a non-saturated model or decreases the level. In 
this section we show that the maximum number of iterations 
required is $5$. (In our {\sf MAGMA} implementation
we count the use of Theorem~\ref{thm:sat} to decrease the level
of a non-saturated model as a further iteration. With this
convention the maximum number of iterations is $6$.)

\begin{Lemma} 
\label{lem:wt29}
  Let $(r,s) \in \W$ be a weight. Then $(r,s)$ dominates one of the 
  weights $w_1, \ldots, w_{29}$ in the following table. (The weights 
  in Lemma~\ref{lem:sevenwts} appear with new numberings. 
  We have marked these weights in bold.)
\[ \hspace{-0.6em} \begin{array}{c|ccccc|ccccc|c||c|ccccc|ccccc|c}
 & \! r_1 \! & \! r_2 \! & \! r_3 \! & \! r_4 \! & \! r_5 \! &
   \! s_1 \! & \! s_2 \! & \! s_3 \! & \! s_4 \! & \! s_5 \! & \la_\nu &
 & \! r_1 \! & \! r_2 \! & \! r_3 \! & \! r_4 \! & \! r_5 \! &
   \! s_1 \! & \! s_2 \! & \! s_3 \! & \! s_4 \! & \! s_5 \! & \la_\nu
\\ \hline
w_{1} & 0 & 0 & 0 & 0 & 0 & \!\! -1 \!\! & 0 & 0 & 0 & 0 & 1  & 
w_{16} & 0 & 1 & 1 & 2 & 2 & 1 & 2 & 2 & 3 & 3 & 7 \\
{\bf{w_{2}}} & 0 & 0 & 0 & 0 & 1 & 0 & 0 & 0 & 0 & 1 & 1  & 
w_{17} & 0 & 1 & 1 & 2 & 2 & 1 & 2 & 2 & 2 & 4 & 6 \\
{\bf{w_{3}}} & 0 & 0 & 1 & 1 & 1 & 1 & 1 & 1 & 1 & 1 & 1  & 
w_{18} & 0 & 1 & 1 & 2 & 2 & 1 & 1 & 2 & 3 & 4 & 7 \\
w_{4} & 0 & 1 & 1 & 1 & 1 & 1 & 1 & 1 & 2 & 2 & 1  & 
{\bf{w_{19}}} & 0 & 1 & 1 & 2 & 3 & 2 & 2 & 3 & 3 & 3 & 6 \\
w_{5} & 0 & 0 & 0 & 1 & 1 & 0 & 0 & 1 & 1 & 1 & 3  & 
{\bf{w_{20}}} & 0 & 1 & 1 & 2 & 3 & 2 & 2 & 2 & 3 & 4 & 7 \\
w_{6} & 0 & 0 & 0 & 1 & 1 & 0 & 0 & 0 & 1 & 2 & 3  & 
w_{21} & 0 & 1 & 1 & 2 & 3 & 1 & 2 & 3 & 3 & 4 & 13 \\
w_{7} & 0 & 0 & 1 & 1 & 1 & 0 & 0 & 1 & 2 & 2 & 3  & 
w_{22} & 0 & 1 & 1 & 2 & 3 & 1 & 2 & 2 & 3 & 5 & 12 \\
w_{8} & 0 & 0 & 1 & 1 & 1 & 0 & 1 & 1 & 1 & 2 & 3  & 
w_{23} & 0 & 1 & 2 & 2 & 3 & 2 & 3 & 3 & 3 & 4 & 9 \\
{\bf{w_{9}}} & 0 & 1 & 1 & 2 & 2 & 2 & 2 & 2 & 2 & 3 & 3  & 
w_{24} & 0 & 1 & 2 & 2 & 3 & 2 & 2 & 3 & 4 & 4 & 9 \\
{\bf{w_{10}}} & 0 & 0 & 1 & 1 & 2 & 1 & 1 & 1 & 2 & 2 & 4  & 
w_{25} & 0 & 1 & 2 & 2 & 3 & 1 & 3 & 3 & 4 & 4 & 10 \\
w_{11} & 0 & 0 & 1 & 1 & 2 & 0 & 0 & 2 & 2 & 3 & 5  & 
w_{26} & 0 & 1 & 2 & 2 & 3 & 1 & 2 & 3 & 4 & 5 & 15 \\
w_{12} & 0 & 0 & 1 & 1 & 2 & 0 & 1 & 2 & 2 & 2 & 8  & 
{\bf{w_{27}}} & 0 & 1 & 2 & 3 & 4 & 3 & 3 & 4 & 4 & 5 & 12 \\
w_{13} & 0 & 0 & 1 & 1 & 2 & 0 & 1 & 1 & 2 & 3 & 8  & 
w_{28} & 0 & 1 & 2 & 3 & 4 & 2 & 3 & 4 & 5 & 5 & 20 \\
w_{14} & 0 & 1 & 1 & 1 & 2 & 1 & 2 & 2 & 2 & 2 & 4  & 
w_{29} & 0 & 1 & 2 & 3 & 4 & 1 & 3 & 4 & 5 & 6 & 22 \\
w_{15} & 0 & 1 & 1 & 1 & 2 & 1 & 1 & 2 & 2 & 3 & 4   
\end{array} \hspace{-0.6em} \]
\normalsize
\end{Lemma}

\begin{Proof}
  We define a {\em standard inequality}
  to be an inequality of the form $r_i + r_j \le s_k + m$ where $1 \le
  i < j \le 5$, $1 \le k \le 5$ and $m$ is a non-negative integer. The
  condition that $(r,s) \in \W$ does not dominate $w_\nu$ is equivalent
  to a list of $\la_\nu$ standard inequalities, at least one of which must
  hold, where $\la_\nu$ is as given in the table. For example, 
  $(r,s) \not\ge w_1$ if and only if $r_1 + r_2 \le s_1$, whereas
  $(r,s) \not\ge w_5$ if and only if $r_1 + r_4 \le s_2$ or 
  $r_4 + r_5 \le s_2 + 1$ or $r_4 + r_5 \le s_5$.
  (We have used the
  conditions $r_1 \le \ldots \le r_5$ and $s_1 \le \ldots \le s_5$ to 
  remove redundant inequalities.)

  We wrote a program using the simplex algorithm to maximise $\sum (2
  r_i -s_i)$ for $(r,s) \in \R^{10}$ subject to $0 \le r_1 \le
  \ldots \le r_5$, $0 \le s_1 \le \ldots \le s_5,$ and a list
  of standard inequalities. Our program starts with the basic feasible 
  solution $(r,s) = (0,0)$. If there is a finite maximum, and it is 
  less than $1$, then (by definition of $\W$) there are no weights 
  satisfying these inequalities. If the maximum is $1$ then we add the 
  constraint $\sum (2 r_i -s_i) =1$. We then use the simplex algorithm 
  to maximise each of the functions $r_i + r_j - s_k$ in turn. 
  In the case of a finite maximum $\alpha$ we obtain an additional 
  standard inequality $r_i + r_j \le s_k + \max( \lfloor \alpha \rfloor, 0)$.
  Then running our original program on the enlarged set of standard
  inequalities we may still be able to show that $\sum (2 r_i -s_i) < 1$.

  After processing the inequalities coming from $w_1, \ldots, w_\nu$
  for $\nu =1, \ldots, 29$ the number of cases remaining were as follows:
  \[ \begin{array}{c} 1, \, 1, \, 1, \, 1, \, 3, \, 5, \, 8, \, 
   13, \, 16, \, 30, \, 31, \, 49, \, 58, \, 47, \, 60, \, 
\qquad \\ \qquad
  64, \, 58, \, 53, \, 45, \, 36, \, 39, \, 34, \, 25, \, 
  15, \, 14, \, 10, \, 3, \, 1, \, 0. \end{array} \]
  The final zero indicates that no cases remain, and this proves
  the lemma.
  The proof of Lemma~\ref{lem:sevenwts} 
  is similar but easier.
\end{Proof}

If $\Phi \in X_5(\OK)$ is non-minimal then by Lemmas~\ref{lem:wts} 
and~\ref{lem:wt29} it has slope at most $14$. This already shows that 
the algorithm in Theorem~\ref{minalgthm}(iii) takes at most $14$ 
iterations. The next lemma improves this bound to $7$ iterations.

\begin{Lemma}
\label{lem:drop}
If the procedure in Theorem~\ref{minalgthm}(ii) returns a saturated model
with the same level then the slope goes down by at least $2$.
\end{Lemma}
\begin{Proof} We revisit the proof of Theorem~\ref{minalgthm}(iii) 
at the end of Section~\ref{sec:wtsl}. If the slope goes down by only 
one then $\Sing \CC_\phi$ spans a hyperplane. If $\Sing \CC_\phi$ is 
a hyperplane then the proof of Theorem~\ref{minalgthm}(ii) at the 
end of Section~\ref{sec:sing} shows that the level is decreased.
Otherwise by Lemma~\ref{lem:sing} we may assume 
$\phi$ takes the form~(\ref{x45}). We then follow the proof of 
Lemma~\ref{lem:w34}(i) with $m=4$. After applying the
transformation suggested there, the second row of $\phi$ has at most
one non-zero entry. This implies that $\Phi$ is non-saturated.
\end{Proof}

The next lemma will be used to show that only $5$ iterations are
required.

\begin{Lemma}
\label{lem:cusp}
Let $\Phi \in X_5(\OK)$ be non-minimal and of slope greater than $10$.
Then replacing $\Phi$ by an $\OK$-equivalent model we may assume
it has weight $w_{29}$ and the coefficient of $x_k$ in $\Phi_{ij}$ is a 
unit for \[ (i,j,k) \in \{ (1,2,1), (1,4,2), (1,5,3), (2,3,2), 
(2,4,3), (2,5,4), (3,4,4), (3,5,5) \}. \]
\end{Lemma}

\begin{Proof}
By Lemma~\ref{lem:wts} we know that $\Phi$ is $\OK$-equivalent to a model
with one of the $29$ weights listed in Lemma~\ref{lem:wt29}. For all but
one of these weights $(r,s)$ we have $\sum_{i=1}^5 (s_i - s_1) \le 10$.
The remaining case is $w_{29}$. If one of the 
coefficients listed is not a unit then $\Phi$ has weight 
$w_\nu$ for some $\nu \in \{1, 5, 13, 26, 16, 21, 8, 12 \}$.
\end{Proof}

We write $[j, \ldots,5]$ for a linear combination of 
$x_j, \ldots, x_5$, and underline in cases where we know the coefficient
is non-zero. Lemma~\ref{lem:cusp} shows that $\Phi \in X_5(\OK)$ has
reduction $\phi \in X_5(k)$ of the form
\[ \begin{pmatrix}
  0  &[\underline{1} ,2,3,4, 5]&[2 ,3,4, 5]&[\underline{2} ,3,4, 5]
     &[\underline{3}, 4, 5 ]\\
     &  0  &[\underline{2} ,3,4, 5]&[\underline{3}, 4, 5 ]&
 [\underline{4}, 5 ]\\
     &     &  0  &[\underline{4}, 5 ]&[\underline{5}]\\
     &     &     &  0  &  0  \\
     &     &     &     &  0  
\end{pmatrix}. \]
Let $\Pf(\phi) = (p_1, \ldots, p_5)$. By considering the partial
derivatives of $p_1,p_2,p_4$ with respect to $x_1, x_2, x_3$
we see that if $P =(x_1: \ldots: x_5) \in \Sing \CC_\phi$ 
then $x_5=0$. Then since $P \in \CC_\phi$ we have $x_4 = x_3 = x_2 = 0$.
So $(1: 0: \ldots : 0)$ is the unique singular point.

Our algorithm applies the transformation
\[[\Diag(1,1,1,\pi^{-1},\pi^{-1}), \Diag(1,\pi,\pi,\pi,\pi)].\]
The result is a model $\Phi$ with weight 
$w_{26} = (0,1,2,2,3;1,2,3,4,5)$ 
whose reduction $\phi$ takes the form
\[ \begin{pmatrix}
  0  &[\underline{1}]& 0 &[\underline{2} ,3,4, 5]
     &[\underline{3}, 4, 5 ]\\
     &  0  &  0  &[\underline{3}, 4, 5 ]&[\underline{4}, 5 ]\\
     &     &  0  &[\underline{4}, 5 ]&[ \underline{5} ]\\
     &     &     &  0  &[ 5 ]\\
     &     &     &     &  0  
\end{pmatrix}. \]
A calculation similar to that above shows that $\Sing \CC_\phi = 
\{ x_3 = x_4 = x_5 = 0 \}$.

Our algorithm applies the transformation
\[[\Diag(\pi,1,1,1,1), \Diag(\pi^{-1},\pi^{-1},1,1,1)]\]
The result is a model $\Phi$ with weight 
$w_{13}=(0,0,1,1,2;0,1,1,2,3)$ 
whose reduction $\phi$ takes the form
\[ \begin{pmatrix}
  0  &[\underline{1}]&0&[\underline{2}]&0\\
     &  0  & [\underline{2}] &[2,\underline{3},4,5]&[ \underline{4},5 ]\\
     &     &  0  &[ \underline{4}, 5 ]&[ \underline{5} ]\\
     &     &     &  0  &[ 5 ]\\
     &     &     &     &  0  
\end{pmatrix}. \] 
A calculation similar to that above shows that 
$\Sing \CC_\phi = \{ x_2 = x_4 = x_5 = 0 \}$. 

The next transformation 
$[\Diag(1,\pi,1,1,1), \Diag(\pi^{-1},1,\pi^{-1},1,1)]$
gives a model with weight $w_{15}=(0,1,1,1,2;1,1,2,2,3)$. 
So after~$3$ iterations the slope is at most~$4$. It follows 
by Lemma~\ref{lem:drop} that at most $5$ iterations are required.

\begin{Example}
The simplest example of a genus one model satisfying the
conditions of Lemma~\ref{lem:cusp} is
\[ \Phi = \begin{pmatrix}
  0  & x_1 &0& x_2 & x_3\\
     &  0  & x_2 & x_3 & x_4 \\
     &     &  0  & x_4 & x_5 \\
     &  -   &     &  0  & 0 \\
     &     &     &     &  0  
\end{pmatrix}. \] 
We find that $\CC_\Phi$ is a rational curve with a cusp, parametrised
by \[ (s:t) \mapsto (-s^5: s^3 t^2 : s^2 t^3 : s t^4 : t^5). \] 
In this case our algorithm takes the maximum of 
exactly $5$ iterations to give a 
non-saturated model. (The first $3$ iterations are already described above.) 
Although the model in this example is singular, 
there are $\pi$-adically close non-singular models that 
are treated in the same way by our algorithm.
\end{Example}

\section{Insoluble models}
\label{sec:insol}

In this section we prove a result converse to the strong 
minimisation theorem. This is analogous to the results for 
models of degrees $n=2,3,4$ proved in \cite[Section 5]{CFS}.
As in Section~\ref{sec:T} we work over a discrete valuation field 
$K$. We write $K^{\sh}$ for the strict Henselisation of $K$. 
(If~$K$ is a $p$-adic field then this is the maximal 
unramified extension.)

\begin{Theorem}
\label{converse}
If $\Phi \in X_5(K)$ is non-singular and $\CC_\Phi(K^\sh) = \emptyset$ then 
the minimal level is at least $1$, and is equal to $1$ if 
$\Char(k) \not= 5$.
\end{Theorem}

As in Section~\ref{sec:numit} 
we write $[j, \ldots,5]$ for a linear combination of 
$x_j, \ldots, x_5$, and underline in cases where we require
the coefficient is non-zero.

\begin{Definition}
\label{def:crit}
A genus one model $\Phi \in X_5(\OK)$ is {\em critical} if it has
reduction mod $\pi$ of the form
\[ \begin{pmatrix}
  0  &[\underline{1},2,3,4,5]&[\underline{2},3,4,5]
&[\underline{3} ,4, 5]
     &[\underline{4}, 5 ]\\
    &  0  &[\underline{3} ,4, 5] & 
 [\underline{4}, 5 ] & [\underline{5}] \\
     &     &  0  &[\underline{5}] & 0 \\
     &     &     &  0  &  0  \\
     &     &     &     &  0  
\end{pmatrix} \]
and $\pi^{-1} \Phi_{35}$, $\pi^{-1} \Phi_{45}$ have reductions mod $\pi$ 
of the form $[\underline{1},2,3,4,5]$, $[\underline{2},3,4,5]$.
\end{Definition}

We show in the next three lemmas that critical models are insoluble,
minimal and of positive level. We then show that every insoluble
model $\Phi \in X_5(K)$ is $K$-equivalent to a critical model.

\begin{Lemma}
Critical models are insoluble over $K$.
\end{Lemma}
\begin{Proof}
Suppose $(x_1, \ldots, x_5) \in K^5$ is a non-zero solution with
$\min \{v(x_i)\} = 0$. By considering the $4 \times 4$ Pfaffians
we successively deduce $\pi \dv x_5$, $\pi \dv x_4$, \ldots, $\pi \dv x_1$.
In particular $\min \{v(x_i)\} > 0$. This is the required contradiction.
\end{Proof}
Since the definition of a critical model is unchanged by an unramified
field extension, it follows immediately that critical models are
insoluble over $K^\sh$.

\begin{Lemma}
Critical models are minimal.
\end{Lemma}
\begin{Proof}
It is easy to see that critical models are saturated. Moreover 
every point on $\CC_\phi = \{ x_3 = x_4 = x_5 = 0 \}$ is singular.
Our algorithm (see Theorem~\ref{minalgthm}) makes the transformation
$[\Diag(\pi,1,1,1,1),\pi^{-1} \Diag(1,1,\pi,\pi,\pi)]$.
This gives an integral model of 
the same level, that is $\OK$-equivalent 
(by a pair of cyclic permutation matrices) to a critical model. 

If $\Phi$ were non-minimal then our algorithm would succeed in reducing
the level. But on the contrary, when given a critical model, our algorithm
endlessly cycles between five $\OK$-equivalence classes.
\end{Proof}

The next lemma describes the possible levels of a critical model.
To treat the cases $\Char(k)=2,3$ 
we need to work with the $a$-invariants defined 
in Section~\ref{g1mod}. Although these are not 
$\SL_5 \times \SL_5$-invariant, if we make our choices of $a_1, b_2, a_3$ 
so as not to introduce any new monomials when we lift to characteristic $0$,
then they will be invariant under all pairs of diagonal matrices. 
It follows by the proof 
of Lemma~\ref{lem:agen} that $a_1, \ldots, a_6$ are isobaric, i.e.
\[ a_i \circ [\Diag(\la_1, \ldots, \la_5), \Diag(\mu_1, \ldots, \mu_5)] = 
(\prod \la_\nu)^{2i} (\prod \mu_\nu)^i a_i. \]

\begin{Lemma}
\label{lem:critlev}
The level of a critical model is at least $1$ and equal to $1$ if 
$\Char(k) \not= 5$.
\end{Lemma}
\begin{Proof}
Applying $[ \Diag(1,\pi^{-1/5},\pi^{-2/5},\pi^{-3/5},\pi^{-4/5}),
  \Diag(\pi^{1/5},\pi^{2/5},\pi^{3/5},\pi^{4/5},\pi) ]$
to a critical model $\Phi$ 
gives a model with coefficients in $\OK[\pi^{1/5}]$. 
It follows by the isobaric property that $\pi^i \dv a_i(\Phi)$ for all $i$.
Hence $\Phi$ has positive level.

The model with coefficients in $\OK[\pi^{1/5}]$ has reduction
\[ \begin{pmatrix}
  0  & \la_1 x_1 & \mu_2 x_2 & - \mu_3 x_3 & -\la_4 x_4 \\
    &  0  & \la_3 x_3 & \mu_4 x_4 & - \mu_5 x_5 \\
     &     &  0  & \la_5 x_5  & \mu_1 x_1 \\
     &     &     &  0  &  \la_2 x_2  \\
     &     &     &     &  0  
\end{pmatrix} \]
for some $\la_1, \ldots, \la_5, \mu_1 \ldots, \mu_5 \in k^\times$. The
invariants of this model are
\begin{align*}
c_4(\la,\mu) 
&= \la^4 + 228 \la^3 \mu + 494 \la^2 \mu^2 - 228 \la \mu^3 + \mu^4, \\
c_6(\la,\mu)
 &= -\la^6 + 522 \la^5 \mu + 10005 \la^4 \mu^2 + 10005 \la^2 \mu^4 
  - 522 \la \mu^5 - \mu^6, 
\end{align*}
and $\Delta(\la,\mu) = \la \mu (\la^2 - 11 \la \mu - \mu^2)^5$,
where $\la = \prod \la_i$ and $\mu = \prod \mu_i$.
Computing a resultant shows that if $\Char(k) \not= 5$ then 
$c_4(\la,\mu)$ and $\Delta(\la,\mu)$ have no common roots.
Therefore the critical model $\Phi$ we started with 
satisfies either $v(c_4(\Phi))=4$ or $v(\Delta(\Phi))=12$.
It follows that $\Phi$ has level at most $1$.
\end{Proof}

\begin{Remark}
The following example of a critical model of level $2$ over $K = \Q_5$
shows that the hypothesis $\Char(k) \not=5$ cannot be removed
from Lemma~\ref{lem:critlev}.
\[ \begin{pmatrix}
    0   &  x_1  &   x_2    & -x_3  &  -x_4 \\
        &    0  &   x_3    & x_4  &  -x_5 \\
        &       &     0    & x_5  & 35 x_1 \\
        &  -    &          &   0  & 5 x_2 \\
        &       &          &      &    0
\end{pmatrix} \]
\end{Remark}

We recall that the minimal level is unchanged by an unramified 
field extension. Replacing $K$ by $K^\sh$ we may assume for the rest of
this section that $K$ is Henselian and its
residue field $k$ is algebraically closed. 
To complete the proof of Theorem~\ref{converse} we show

\begin{Theorem}
\label{showcrit}
If $\Phi \in X_5(\OK)$ is minimal
and $\CC_\Phi(K) = \emptyset$ then $\Phi$ is $\OK$-equivalent to a
critical model.
\end{Theorem}

We start the proof of Theorem~\ref{showcrit} with the following lemma.

\begin{Lemma}
\label{minlem}
If $\Phi \in X_5(\OK)$ is minimal then its reduction 
$\phi \in X_5(k)$ has the following properties.
\begin{enumerate}
\item The $4 \times 4$ Pfaffians of $\phi$ are linearly independent.
\item The subscheme $\CC_\phi \subset \PP^4$ does not contain a plane.
\item The entries of $\phi$ span the space of linear forms on $\PP^4$.
\end{enumerate}
\end{Lemma}
\begin{Proof}
(i) This follows by Theorem~\ref{thm:sat} and Lemma~\ref{lem:sat}(i). \\
(ii) Suppose $\CC_\phi$ contains the plane $\{ x_4 = x_5 = 0 \}$.
By Lemma~\ref{lem:pfzero} we may assume the reduction mod
$\pi$ of $[I_5,\Diag(1,1,1,\pi,\pi)]\Phi$ takes one of the
two forms given in the lemma.
We decrease the level by applying either 
$[\Diag(\pi,1,1,1,1),\pi^{-1}I_5]$ or $[\Diag(1,1,1,\pi^{-1},\pi^{-1}),B]$ 
where $B$ is chosen to preserve integrality. \\
(iii) This is clear, as we could otherwise decrease the level by 
dividing one of the co-ordinates by $\pi$. 
\end{Proof}

\begin{Lemma} 
\label{threecases}
Let $\phi \in X_5(k)$ be a genus one model satisfying
the conclusions of Lemma~\ref{minlem}. Suppose that every point on
$\CC_\phi$ is singular. Then $\phi$ is $k$-equivalent to 
\begin{equation*}
 \begin{pmatrix}
    0   &  0  &   x_1    & x_3  &  x_4 \\
        &   0  &   x_2    & x_4  &  x_5 \\
        &       &     0    & x_5  & 0 \\
        &  -    &          &   0  & 0 \\
        &       &          &      &    0
\end{pmatrix} 
\text{ or }
\begin{pmatrix}
    0   &   x_1   &   0    & x_3  &  x_4 \\
        &   0  &   x_2    & x_4  &  x_5 \\
        &       &     0    & x_5  & 0 \\
        &  -    &          &   0  & 0 \\
        &       &          &      &    0
\end{pmatrix} 
\text{ or }
\begin{pmatrix}
    0   &  x_1  &   x_2    & x_3  &  x_4 \\
        &   0  &   x_3    & x_4  &  x_5 \\
        &       &     0    & x_5  & 0 \\
        &  -    &          &   0  & 0 \\
        &       &          &      &    0
\end{pmatrix}. 
\end{equation*}
\end{Lemma}

Our proof of Lemma~\ref{threecases} uses the following 
classification of degenerations of the twisted cubic.
(Only the last sentence of the statement is needed.)

\begin{Lemma}
\label{lem:tc}
Let $\psi$ be a $3 \times 2$ matrix of linear forms in $R = k[x_1,\ldots,x_4]$.
Suppose the $2 \times 2$ minors of $\psi$ are linearly independent and
no linear combination of them has rank $1$. Then $\psi$ is 
$\GL_2 \times \GL_3 \times \GL_4$-equivalent to one of the following:
\begin{equation}
\label{tcdegen}
\begin{pmatrix}
    x_1 & x_2 \\ x_2 & x_3 \\ x_3 & x_4 
\end{pmatrix}, \quad \begin{pmatrix}
    x_1 & x_2 \\ x_2 & x_3 \\ x_4 & 0 
\end{pmatrix}, \quad \begin{pmatrix}
    x_1 & x_2 \\ 0 & x_3 \\ x_4 & 0 
\end{pmatrix}, \quad \begin{pmatrix}
    x_1 & 0 \\ x_2 & x_2 \\ 0 & x_3 
\end{pmatrix}. 
\end{equation}
In particular the locus of smooth points on 
$\Gamma = \{ \rank \psi \le 1\} \subset \PP^3$ spans $\PP^3$.
\end{Lemma}
\begin{Proof} We may realise $\Gamma$ as the intersection of the image
of the Segre embedding $\PP^1 \times \PP^2 \to \PP^5$ with a
linear subspace $\PP^3$. So every component of $\Gamma$ has dimension
at least $1$. If every component has dimension $1$ then by the
Buchsbaum-Eisenbud acyclicity criterion there is a minimal free
resolution 
\begin{equation}
\label{resol}
 0 \ra R(-3)^2 \stackrel{\psi}{\ra} R(-2)^3 \stackrel{M}{\ra} R 
\end{equation}
where 
$M$ is the vector of 
$2 \times 2$ minors of $\psi$. If in addition $\dim T_P \, \Gamma =1$
for every $P \in \Gamma$ then by an argument using
Serre's criterion (see \cite[Section 18.3]{E})
the ideal in $R$ generated by the $2 \times 2$ minors of $\psi$ 
is a prime ideal.
By~(\ref{resol}) the Hilbert polynomial is
\[ h(t) = \binom{t+3}{3} -3 \binom{t+1}{3} + 2 \binom{t}{3} = 3t+1. \]
Therefore $\Gamma$ is a twisted cubic and $\psi$ is equivalent to 
the first of the matrices in~(\ref{tcdegen}).

In all other cases there exists $P \in \Gamma$ with $\dim T_P \, \Gamma > 1$.
First suppose $\rank \psi(P) = 1$. Moving $P$ to $(1:0:0:0)$ we may
suppose 
\[ \psi = \begin{pmatrix} x_1 & \alpha \\ \delta & \beta \\ \gamma & 0
\end{pmatrix} \]
where $\alpha,\beta,\gamma,\delta$ are linear forms in $x_2,x_3,x_4$.
Our hypotheses on the $2 \times 2$ minors ensure that $\alpha,\beta,\gamma$
are linearly independent; say they are $x_2,x_3,x_4$.
By row and column operations (and a substitution for $x_1$)
we may assume $\delta$ is a multiple of $x_2$. This gives
the second and third cases in~(\ref{tcdegen}).

Now suppose $\rank \psi(P) = 0$. Let $Q \in \Gamma$ be any other point. 
If $\rank \psi(Q) = 0$ then the $2 \times 2$ minors are binary quadratic
forms, and so some linear combination has rank $1$. Therefore
$\rank \psi(Q) = 1$. If $\dim T_Q \, \Gamma > 1$ then our earlier analysis
applies (and in fact gives a contradiction). Otherwise we may assume
\[ \psi = \begin{pmatrix} x_1 & 0 \\ \alpha & x_2 \\ \beta & x_3
\end{pmatrix} \]
where $\alpha,\beta$ are linear forms in $x_2,x_3$. 
(The zero in the top right has been cleared by row operations.) Since 
$\alpha x_3 - \beta x_2$ is a rank~$2$ quadratic form in $x_2,x_3$
we can make a change co-ordinates so that 
$\Gamma = \{ x_1 x_2 = x_1 x_3 = x_2 x_3 = 0 \}$. Then
$\psi$ is equivalent to the last of the matrices in~(\ref{tcdegen}). 

For the final statement, we note that the $4$ cases correspond 
geometrically to (i) a twisted cubic, (ii) a conic and a line,
(iii) three non-concurrent lines, and (iv) three concurrent lines.
In each case $\Gamma$ spans $\PP^3$ and the only singular points
are the points where the components meet.
\end{Proof}

\begin{ProofOf}{Lemma~\ref{threecases}}
Let $P \in \CC_\phi$ be a singular point. Moving $P$ to $(1:0:0:0:0)$
we may assume $\phi$ takes the form
\[ \begin{pmatrix}
    0   &  x_1  &   \ell_2    & \alpha_1  &  \beta_1 \\
        &   0  &   \ell_3    & \alpha_2  &  \beta_2 \\
        &       &     0    & \alpha_3  & \beta_3 \\
        &  -    &          &   0  & 0 \\
        &       &          &      &    0
\end{pmatrix} \]
where $\ell_i, \alpha_i, \beta_i$ are linear forms in $x_2, \ldots, x_5$.
Let $\psi$ be the top right $3 \times 2$ submatrix and let
$\Gamma \subset \PP^3$ be the curve defined by its $2 \times 2$ minors. 
Since the $2 \times 2$
minors of $\psi$ are a subset of the $4 \times 4$ Pfaffians of $\phi$,
they are linearly independent. In particular $\alpha_3$ and 
$\beta_3$ cannot both vanish identically. Without loss of 
generality $\alpha_3$ is non-zero. 

Suppose no linear combination of the $2 \times 2$ minors of $\psi$ has 
rank $1$. Then by Lemma~\ref{lem:tc} there is a smooth point 
$Q =(x_2:x_3:x_4:x_5)$ on $\Gamma$ with $\alpha_3(Q) \not= 0$. 
Solving for $x_1$ gives a smooth point $(x_1:x_2:\ldots:x_5)$
on $\CC_\phi$. This is a contradiction. 
Therefore some linear
combination of the $2 \times 2$ minors of $\psi$ has 
rank $1$. It is then easy to see that $\phi$ is $k$-equivalent to a 
model of the form~(\ref{x45}).

By properties (i) and (ii), $\eta, x_5$ are linearly independent
and $\gamma,\delta,x_5$ are linearly independent. However if 
$\eta,\gamma,\delta,x_5$ were linearly independent then taking them
to be $x_2, \ldots, x_5$ would give that $(0:1:0:0:0)$ is
a smooth point on $\CC_\phi$. By row and column operations we
may therefore suppose $\eta = \delta$ ($= x_4$ say). 

By property (ii), $\beta, x_4, x_5$ are linearly 
independent and $\gamma, x_4, x_5$ are linearly independent.
By row and column operations (and substitutions for the $x_i$)
we may suppose $\beta = x_3$ and $\gamma = x_2$ or $x_3$.
If $\gamma = x_2$ then by further row and column operations 
(and substitutions for the $x_i$) we may suppose $\alpha$ is a multiple
of $x_1$. The lemma now follows using property (iii).
\end{ProofOf}

\begin{ProofOf}{Theorem~\ref{showcrit}}
Since $K$ is Henselian any smooth point on $\CC_\phi$ lifts
to a $K$-point on $\CC_\Phi$. So 
we may assume $\phi$
takes one of the three forms in Lemma~\ref{threecases}. In the first
two cases $\phi$ defines a pair of concurrent lines 
with multiplicities~$2$ and~$3$. (These cases may be distinguished
by the dimension of the tangent space at the point of intersection).
In the third case it defines a line with multiplicity~$5$.

We apply the transformation $[\Diag(1,1,1,1,\pi^{-1}),\Diag(1,1,1,\pi,\pi)]$.
This gives an integral model of the same level.
So the reduction must again be $k$-equivalent to one of the
three models in Lemma~\ref{threecases}.
We tidy up by an $\OK$-equivalence 
that cyclically permutes the rows 
and columns, and makes substitutions for $x_4$ and $x_5$.
The reduction $\phi \in X_5(k)$ now takes the form
\[ \begin{pmatrix}
    0   &  x_4  &  x_5    & \alpha  & \beta \\
        &   0  &   0    & x_1  &  x_3 \\
        &       &     0    &  x_2  & 0 \\
        &  -    &          &   0  & 0 \\
        &       &          &      &    0
\end{pmatrix} \text{ or } \begin{pmatrix}
    0   &  x_4  &  x_5    & \alpha  & \beta \\
        &   0  &   x_1    & 0  &  x_3 \\
        &       &     0    &  x_2  & 0 \\
        &  -    &          &   0  & 0 \\
        &       &          &      &    0
\end{pmatrix} \text{ or } 
\begin{pmatrix}
    0   &  x_4  &  x_5    & \alpha  & \beta \\
        &   0  &   x_1    & x_2  &  x_3 \\
        &       &     0    &  x_3  & 0 \\
        &  -    &          &   0  & 0 \\
        &       &          &      &    0
\end{pmatrix} \]
where $\alpha$ and $\beta$ are linear forms in $x_1,x_2,x_3$.

In the first case $(0:0:0:1:0)$ is a point with tangent space of
dimension $3$ and $\CC_\phi$ contains points not on the line 
$\{x_1 = x_2 = x_3 = 0\}$. So the transformation has moved us to the
second case.

In second case we obtain a contradiction as follows.
If $\alpha = x_1 + \la x_2 + \mu x_3$ then
adding $\mu$ times the fifth row/column to the third
row/column, and making substitutions for $x_1$ and $x_5$ we may assume
$\mu = 0$. Then $(0:0:1:0:0)$ is a smooth point on $\CC_\phi$.
Likewise if $\beta = x_1 + \la x_2 + \mu x_3$ then subtracting
$\la$ times the fourth row/column from the second row/column
and making substitutions for $x_1$ and $x_4$ we may assume 
$\la = 0$. Then $(0:1:0:0:0)$ is a smooth point on $\CC_\phi$.
We are forced to the conclusion that neither $\alpha$ nor
$\beta$ involves $x_1$. But then $\CC_\phi$ contains the plane
$\{x_2 = x_3=0\}$ and by Lemma~\ref{minlem} this contradicts that
$\Phi$ is minimal.

In the third case we show that if the transformation above
brings us back to the third case, then the original model 
is critical. If $\beta = x_1 + \la x_2 + \mu x_3$
then adding $\la$ times the fourth row/column to the third row/column, 
and making substitutions for $x_1$ and $x_5$ we may 
assume $\la =0$. Then $\CC_\phi$ contains the lines
$\{x_1 = x_2 = x_3 = 0\}$ and $\{x_1 = x_3 = x_5 = 0\}$. 
So if the transformation returns us to third case then $\beta$ cannot 
involve $x_1$. Since $\CC_\phi$ does not contain a plane, and the 
$4 \times 4$ Pfaffians of $\phi$ are linearly independent, $\alpha$ 
must involve $x_1$ and $\beta$ must involve $x_2$. 
It follows by Definition~\ref{def:crit} that the original model is
$\OK$-equivalent to a critical model.
\end{ProofOf}

\section{Reduction} 
\label{sec:R}

Let $C \subset \PP^{4}$ be a genus one normal curve of degree
$5$ defined over $\Q$. We may represent it by a non-singular 
genus one model $\Phi \in X_5(\Z)$. Running the algorithm
in Section~\ref{sec:minalg} locally at $p$, for all primes $p$ dividing 
the discriminant $\Delta(\Phi)$, we obtain a $\Q$-equivalent model
(still with coefficients in $\Z$)
whose discriminant is minimal in absolute value. If 
$C$ is everywhere locally soluble then this discriminant
is the minimal discriminant of $E = \Jac(C)$. 
It remains to make a $\GL_5(\Z)$ change of co-ordinates on $\PP^4$ 
so that (after running the LLL algorithm on the space of $5$ quadrics 
defining the curve) the coefficients (and not just the 
invariants) are small. The general method, described 
in~\cite[Section 6]{CFS}, 
is to run the LLL algorithm on the Gram matrix for the (unique) 
Heisenberg invariant inner product. In this section we outline how to 
compute this inner product in the case $n=5$.

We recall that the Heisenberg group is the subgroup of $\SL_5(\C)$
consisting of matrices $M_T$ that describe the action of $T \in E[5]$
on $C \subset \PP^4$ by translation. For $T \not= 0_E$ we call the $5$ 
points in $\PP^4$ fixed by $M_T$ a {\em syzygetic $5$-tuple}.
It may be shown (for example by adapting the proof of 
\cite[Proposition 4.1]{g1hessians} or using that $H^1(\R,E[5])$ is trivial) 
that $\Phi$ is $\SL_5(\R) \times \SL_5(\R)$-equivalent
to a model in Hesse form:
\begin{equation}
\label{hesseform}
 \begin{pmatrix}
 0  & a x_0 & b x_1 & -b x_2 & -a x_3 \\
 & 0 & a x_2 &  b x_3 & -b x_4 \\
 && 0 & a x_4 & b x_0 \\
 &-&& 0 & a x_1 \\
 &&&& 0 
\end{pmatrix}. 
\end{equation}
The invariants of this model are
\begin{align*}
c_4 &= 
a^{20} + 228 a^{15} b^5 + 494 a^{10} b^{10} - 228 a^5 b^{15} + b^{20}, \\
c_6 &=
-a^{30} + 522 a^{25} b^5 + 10005 a^{20} b^{10} + 10005 a^{10} b^{20} 
- 522 a^5 b^{25} - b^{30},
\end{align*}
and $\Delta = D^5$ where $D = ab(a^{10} - 11 a^5 b^5 - b^{10})$.
For a model in Hesse form the Heisenberg group is generated by
$\Diag(1, \zeta, \ldots, \zeta^4)$, where $\zeta$ is
a primitive $5$th root of unity, and a cyclic permutation matrix.
Since these matrices are unitary, the Heisenberg invariant inner 
product is the standard inner product on $\R^5$.

The Hessian, introduced in \cite{g1hessians}, 
is an $\SL_5 \times \SL_5$-equivariant polynomial map
$H: X_5 \to X_5$ with the property that the Hessian of~(\ref{hesseform})
is of the same form with $a$ and $b$ replaced by $-\partial D/ \partial b$
and $\partial D/\partial a$.

\begin{Theorem}
\label{thm:syz}
Let $\Phi \in X_5(\C)$ be a non-singular genus one model 
with invariants $c_4$ and~$c_6$. 
Let $A$ be the $3 \times 5$ matrix of quadrics 
such that $\la \Phi + \mu H(\Phi)$ has $4 \times 4$ Pfaffians 
\[\{ \la^2 A_{1i} + \la \mu A_{2i} + \mu^2 A_{3i} : i= 1,\ldots, 5\}\]
Then $\X = \{ \rank A \le 1 \} \subset \PP^4$ 
consists of $30$ points and the syzygetic $5$-tuples for
$\CC_\Phi$ are the fibres of the map $\alpha : \X \to \PP^2$ given
by the first (or indeed any) column of~$A$. The image of $\alpha$ 
is the set of $6$ points $(x:y:z) \in \PP^2$ satisfying
\begin{equation}
\label{mat1}
\rank \begin{pmatrix}
0 & 5 x & y & 6 c_4 x + z \\
x & y & 6 c_4 x - z & 8 c_6 x \\
y & -z & 8 c_6 x & 9 c_4^2 x
\end{pmatrix} \le 2. 
\end{equation}
\end{Theorem}
\begin{Proof} It suffices to prove this for $\Phi$ in Hesse form.
Then $\X$ is defined by
\begin{equation}
\label{mat35}
\rank \begin{pmatrix}
x_0^2 & x_1^2 & x_2^2 & x_3^2 & x_4^2 \\
x_1 x_4 & x_0 x_2 & x_1 x_3 & x_2 x_4 & x_0 x_3 \\
x_2 x_3 & x_3 x_4 & x_0 x_4 & x_0 x_1 & x_1 x_2 
\end{pmatrix}  \le 1 
\end{equation}
and by \cite[Proposition 1]{BHM} is a set of 30 points. 
Evaluating the columns of~(\ref{mat35}) at these points 
we obtain $(1:0:0)$ and $(1:\zeta^{i}: \zeta^{-i})$ 
for $i=0,\ldots,4$. 
These are the points $(\xi:\eta:\nu) \in \PP^2$ satisfying
\begin{equation}
\label{mat2}
\rank \begin{pmatrix} 
\xi & \eta & \nu & 0 \\ \nu & \xi & 0 & -\eta \\ 
0 & 0 & \eta & \nu \end{pmatrix} \le 2. 
\end{equation}

The remaining statements follow by direct calculation.
In particular our description~(\ref{mat1}) of the image of 
$\alpha$ is checked by making the substitution
\[ \begin{pmatrix} x \\ y \\ z \end{pmatrix}
= \begin{pmatrix}
a b & b^2 & -a^2 \\ 
-a \frac{\partial D}{\partial a}+b \frac{\partial D}{\partial b} 
& -2 b \frac{\partial D}{\partial a} & -2 a \frac{\partial D}{\partial b} \\ 
-\frac{\partial D}{\partial b} \frac{\partial D}{\partial a} & 
(\frac{\partial D}{\partial a})^2 & -(\frac{\partial D}{\partial b})^2
\end{pmatrix}
 \begin{pmatrix} \xi \\ \eta \\ \nu \end{pmatrix}. \]
We note that this change of co-ordinates, and the matrix
relating the $3 \times 3$ minors of~(\ref{mat1}) and~(\ref{mat2}), 
each have determinant a constant times a power of $D$.
\end{Proof}

After computing the Hessian exactly (using the algorithm in 
\cite[Section 11]{g1hessians}) we use Theorem~\ref{thm:syz} to 
compute the syzygetic $5$-tuples numerically. We then compute
a Gram matrix for the Heisenberg invariant inner product as follows.

\begin{Proposition}
Let $C \subset \PP^4$ a genus one normal curve defined over $\R$.
\begin{enumerate}
\item Exactly two of the syzygetic $5$-tuples for $C$ are defined
over $\R$, say
\begin{align*}
Y & = \{ y_i y_j = 0 : i < j \} \subset \PP^4, \\
Z & = \{ z_i z_j = 0 : i < j \} \subset \PP^4, 
\end{align*}
where $y_0, \ldots, y_4$ and $z_0, \ldots, z_4$ are linear forms
in $\C[x_0, \ldots, x_4]$.
\item One of the $5$-tuples in (i) has $5$ real points and
the other has $1$ real point. We may therefore arrange that
$y_0, \ldots, y_4$ and $z_0$ have real coefficients and that
the pairs $z_1, z_4$ and $z_2, z_3$ are complex conjugates.
\item 
The Heisenberg invariant quadratic form spans 
the $1$-dimensional real vector space 
\[\langle y_0^2, \ldots, y_4^2 \rangle \cap \langle z_0^2, z_1 z_4,
z_2 z_3 \rangle.\]
\end{enumerate}
\end{Proposition}
\begin{Proof}
For $C$ in Hesse form we may take $y_i = x_i$ and 
$z_i = \sum_{j=0}^4 \zeta^{ij} x_j$.
In this case the Heisenberg invariant quadratic form  
is $x_0^2 + \ldots + x_4^2$.
\end{Proof}

\section{Examples}
\label{sec:E}

Wuthrich \cite{W} constructed an element of order $5$
in the Tate-Shafarevich group of the elliptic curve $E/\Q$ with
Weierstrass equation
\[ y^2 + x y + y = x^3 + x^2 - 3146 x + 39049. \]
His example (see also \cite[Section 9]{g1inv}) 
is defined by the $4 \times 4$ Pfaffians of
\small
\[ \begin{pmatrix}
0  & 310 x_1 + 3 x_2 + 162 x_5 &  -34 x_1 - 5 x_2 - 14 x_5  
& 10 x_1 + 28 x_4 + 16 x_5 &  80 x_1 - 32 x_4 \\
&  0 &  6 x_1 + 3 x_2 + 2 x_5 &  -6 x_1 + 7 x_3 - 4 x_4 &  -14 x_2 - 8 x_3 \\
&  &  0 &  -x_3  & 2 x_2 \\
& - &  &  0  & -4 x_1 \\
&  &  &  &  0
\end{pmatrix}. \]
\normalsize
This model has discriminant $2^{132} \Delta_E$ where $\Delta_E$ is the
minimal discriminant of $E$. In other words, the model is minimal
at all primes except $p=2$, where the level is~$11$.
Minimisation and reduction suggest the change of co-ordinates
\small
\[ \begin{pmatrix} x_1 \\ x_2 \\ x_3 \\ x_4 \\ x_5 \end{pmatrix}
\leftarrow
\begin{pmatrix}
  0 &  4  & -8  &  4  &  8 \\
  0 &   0 &   0 &   0 &  16 \\
  0 &  -4  &  4  &  0  & 12 \\
  4 &   5 & -15  &  2  &  7\\
  4 & -12 &  20 & -12  & -8
\end{pmatrix} \begin{pmatrix} x_1 \\ x_2 \\ x_3 \\ x_4 \\ x_5 \end{pmatrix} \]
\normalsize
so that Wuthrich's example simplifies to 
\small
\[ \Phi = \begin{pmatrix}
0 &  x_2 + x_5 &  -x_5  & -x_1 + x_2 &  x_4 \\
 &  0 &  x_2 - x_3 + x_4 &  x_1 + x_2 + x_3 - x_4 - x_5  & x_1 - x_2 - x_3 - x_4 - x_5 \\
 &  &  0 &  x_1 - x_2 + 2 x_3 - x_4 - x_5 &  -x_2 - x_4 + x_5 \\
 & - &  &  0 &  -x_3 - x_4 - 2 x_5 \\
 &  &  &  &  0
\end{pmatrix}. \]
\normalsize
Our {\sf MAGMA} function {\tt DoubleGenusOneModel}, 
described in \cite{invenqI},
computes a genus one model $\Phi'$ that represents twice the class
of $\Phi$ in the $5$-Selmer group. This model has entries
\tiny
\begin{align*}
\Phi'_{12} &= 3534132778 x_1 + 3583651940 x_2 - 881947110 x_3 - 323014538 x_4 + 3395115339 x_5, \\ 
\Phi'_{13} &= 5079379222 x_1 - 2965539950 x_2 + 11022202860 x_3 + 12821590868 x_4 + 640276471 x_5, \\ 
\Phi'_{14} &= -10098238458 x_1 - 1274966110 x_2 - 7873816170 x_3 - 3456923272 x_4 - 62353929 x_5, \\ 
\Phi'_{15} &= -12929747724 x_1 - 6790511810 x_2 - 11113305270 x_3 - 15161763156 x_4 + 3241937033 x_5, \\ 
\Phi'_{23} &= -3381247332 x_1 + 3810679160 x_2 + 5919634530 x_3 + 75326852 x_4 - 1245085426 x_5, \\ 
\Phi'_{24} &= -3572860258 x_1 - 5569480730 x_2 - 953739600 x_3 - 2138046812 x_4 - 858145244 x_5, \\ 
\Phi'_{25} &= -4674149266 x_1 - 943631490 x_2 - 6754488160 x_3 + 751535046 x_4 + 117685567 x_5, \\ 
\Phi'_{34} &= -1851228934 x_1 + 5238146110 x_2 - 165588410 x_3 - 2070411506 x_4 + 678105748 x_5, \\ 
\Phi'_{35} &= -6992835070 x_1 - 3744630360 x_2 + 3130208220 x_3 - 4523781310 x_4 + 433739425 x_5, \\ 
\Phi'_{45} &= 780078472 x_1 + 2039763820 x_2 - 450062790 x_3 - 7105731722 x_4 + 1625466111 x_5.
\end{align*}
\normalsize
The discriminant of $\Phi'$ is $\Delta_E^{49}$. In particular this model
is non-minimal at all bad primes of $E$. Minimisation and
reduction suggest the change of co-ordinates
\small
\[ \begin{pmatrix} x_1 \\ x_2 \\ x_3 \\ x_4 \\ x_5 \end{pmatrix}
\leftarrow
\begin{pmatrix}
  92 & -36 &-153 & 129 &-131 \\
 -54 &  84 &   5 &-206 & 139 \\
 -63 &-174 & -60 & -79 &  53 \\
-111 & 106 & 206 &-115 &-162 \\
 314 &-466 & 158 &-328 & -12
\end{pmatrix} 
\begin{pmatrix} x_1 \\ x_2 \\ x_3 \\ x_4 \\ x_5 \end{pmatrix} \]
\normalsize
so that $\Phi'$ simplifies to 
\small
\[ \begin{pmatrix}
0  & -x_4 + x_5 &  x_3 - x_4 + x_5 &  x_2 - x_5 &  x_1 - x_2 + x_3 - x_4 - 2 x_5 \\
 &  0 &  x_1 + x_5  & -x_2 - x_3 &  -x_2 + x_5 \\
 &  &  0 &  x_4 &  -x_1 \\
 & - &  &  0 &  x_1 + x_4 - x_5 \\
 &  &  &  &  0
\end{pmatrix}. \]


\normalsize

\begin{thebibliography}{MM}

\frenchspacing
\renewcommand{\baselinestretch}{1}

\bibitem[ARVT]{ARVT}
M. Artin, F. Rodriguez-Villegas and J. Tate, 
On the Jacobians of plane cubics,
{\em Adv. Math.} {\bf{198}} (2005), no. 1, 366--382.

\bibitem[BCP]{magma}
W. Bosma, J. Cannon and C. Playoust, 
The {\sf MAGMA} algebra system I: The user language, 
{\em J. Symb. Comb.} {\bf{24}},  (1997) 235--265. 
\,\, \url{http://magma.maths.usyd.edu.au/magma/}

\bibitem[BHM]{BHM}
W. Barth, K. Hulek and R. Moore,
Shioda's modular surface $S(5)$ and the Horrocks-Mumford bundle,
{\em Vector bundles on algebraic varieties} (Bombay, 1984), 35--106,
Tata Inst. Fund. Res. Stud. Math., {\bf{11}}, 
Tata Inst. Fund. Res., Bombay, 1987.

\bibitem[CFS]{CFS}
J.E. Cremona, T.A. Fisher and M. Stoll, 
Minimisation and reduction of 2-, 3- and 4-coverings of elliptic curves, 
{\em Algebra \& Number Theory} {\bf{4}} (2010), no. 6, 763-820.

\bibitem[E]{E}
D. Eisenbud, 
{\em Commutative algebra, with a view toward algebraic geometry},
Graduate Texts in Mathematics {\bf{150}}, Springer-Verlag, New York, 1995.

\bibitem[F1]{g1inv}
T.A. Fisher,
The invariants of a genus one curve, 
{\em Proc. Lond. Math. Soc.} (3) {\bf 97} (2008) 753-782. 

\bibitem[F2]{g1hessians}
T.A. Fisher,
The Hessian of a genus one curve, to appear in {\em Proc. Lond. Math. Soc}.

\bibitem[F3]{invenqI}
T.A. Fisher,
{\em Invariant theory for the elliptic normal quintic, I. Twists of X(5)},
preprint, \url{arXiv:1110.3520v1}

\bibitem[K]{Kraus}
A. Kraus, 
Quelques remarques \`a propos des invariants $c\sb 4,\;c\sb 6$ et 
$\Delta$ d'une courbe elliptique. 
{\em Acta Arith.} {\bf{54}} (1989), no. 1, 75--80.

\bibitem[W]{W}
C. Wuthrich, 
Une quintique de genre 1 qui contredit le principe de Hasse,
{\em Enseign. Math.} (2) {\bf{47}} (2001), no. 1-2, 161--172.

\end{thebibliography}
\end{document}